# EXACT ASYMPTOTICS FOR FLUID QUEUES FED BY MULTIPLE HEAVY-TAILED ON–OFF FLOWS

By Bert Zwart, Sem Borst and Michel Mandjes

*Eindhoven University of Technology and CWI*

We consider a fluid queue fed by multiple On–Off flows with heavy-tailed (regularly varying) On periods. Under fairly mild assumptions, we prove that the workload distribution is asymptotically equivalent to that in a reduced system. The reduced system consists of a "dominant" subset of the flows, with the original service rate subtracted by the mean rate of the other flows. We describe how a dominant set may be determined from a simple knapsack formulation.

The dominant set consists of a "minimally critical" set of On–Off flows with regularly varying On periods. In case the dominant set contains just a single On–Off flow, the exact asymptotics for the reduced system follow from known results. For the case of several On–Off flows, we exploit a powerful intuitive argument to obtain the exact asymptotics. Combined with the reduced-load equivalence, the results for the reduced system provide a characterization of the tail of the workload distribution for a wide range of traffic scenarios.

**1. Introduction.** Over the past few decades, fluid models have gained strong ground as a versatile approach for analyzing burst-scale traffic behavior. The basic model is that of several On–Off sources, each alternating between activity phases (commonly referred to as bursts) and silence periods. When active, each source generates traffic at some constant rate.

Classical papers of Anick, Mitra and Sondhi [2] and Kosten [24] considered a queue fed by the superposition of several homogeneous On–Off sources with exponentially distributed activity and silence periods. Subsequent work extended the model in various directions, such as heterogeneous source characteristics, several source states to account for various activity levels, or activity periods with a general Markovian structure; see, for instance, [25, 38]. Under traditional statistical assumptions, it turns out that the tail of the backlog distribution typically exhibits exponential decay.

---

In recent years, empirical findings have triggered a strong interest in fluid models with non-Markovian activity periods. Extensive measurements indicate that bursty traffic behavior may extend over a wide range of time scales, manifesting itself in long-range dependence and self-similarity; see [26, 33]. The occurrence of these phenomena is commonly attributed to extreme variability and heavy-tailed characteristics in the underlying activity patterns (connection times, file sizes, scene lengths); see [5, 13, 39]. Fluid queues with heavy-tailed activity periods provide a natural paradigm for capturing these characteristics. We refer to [10] for a survey.

Although the presence of heavy-tailed traffic characteristics is widely acknowledged, the practical implications for network performance and traffic engineering remain to be fully resolved. Analytical studies show potentially dramatic performance repercussions for infinite buffers. For moderate buffer sizes, though, the impact of heavy-tailed traffic characteristics is not as pronounced; see [18, 20, 30, 37]. For larger buffer sizes, flow control mechanisms play a critical role in preventing badly behaved traffic from overwhelming the buffer content; see [3]. However, the amount of backlogged traffic at the user, and thus the end-to-end quality-of-service, may still be significantly affected by heavy-tailed activity patterns.

The effect of heavy-tailed traffic characteristics on buffer behavior also crucially depends on the relative amount of heavy-tailed traffic, in particular whether or not activity of heavy-tailed flows alone can cause the buffer to fill. Asymptotic bounds in [15] indeed show a sharp dichotomy in the qualitative behavior of the workload, depending on whether the mean rate of the light-tailed flows plus the peak rate of the heavy-tailed flows exceeds the link rate or not. In case the link rate is larger, the workload distribution has light-tailed characteristics, whereas the link rate being smaller results in heavy-tailed characteristics. The exact asymptotics for the former case were recently obtained in [7]. For the latter case, the bounds of [15] indicate that one can usually identify a "dominant" set, which is a minimal set of flows that can cause a positive drift in the buffer. As far as bounds are concerned, all other flows can essentially be accounted for by subtracting their aggregate mean rate from the link rate. Somewhat related notions are developed in [27] in the setting of $M/G/\infty$ input with heterogeneous sessions. Exact results, however, have remained elusive for all but a few special cases. Results of Agrawal, Makowski and Nain [1] show that the dominance principle described above in fact extends to the exact asymptotics in the case of a *single* dominant flow. This may be expressed in terms of a "reduced-load equivalence," implying that the workload is asymptotically equivalent to that in a reduced system. The reduced system consists only of the dominant flow, with the link rate subtracted by the aggregate mean rate of all other flows. This extends results of Boxma [9], Jelenković and Lazar [22] and Rolski, Schlegel and Schmidt [36] for multiplexing a single (intermediately) regularly varying flow with several exponential flows. Related results are derived in [22, 35] in the context of $M/G/\infty$ input. Like the reduced-load equivalence, however, all these results rely on the assumption that a single active flow is sufficient for a positive drift in the buffer.

In the present paper we determine the exact asymptotics for the case where several On–Off flows must be active for the buffer to fill (under the assumption that the distribution of the On periods is regularly varying [6]). From a practical perspective, this case appears particularly relevant, as the peak rate of a single flow is usually substantially smaller than the link rate. However, the rather subtle interaction of several flows that is involved in filling the buffer drastically complicates the analysis, reflecting the sharp demarcation in known results described above. We start with extending the reduced-load equivalence to the case of a reduced system consisting of several flows, using sample-path arguments. We then build on a qualitative understanding of the large-deviations behavior to obtain the exact asymptotics for the reduced system. This part of the analysis is related to recent work of Resnick and Samorodnitsky [35] on fluid queues with $M/G/\infty$ input.

The remainder of the paper is organized as follows. In Section 2, we present a detailed model description. In Section 3, we give a broad overview of the main results of the paper, and describe how the dominant set may be determined from a simple knapsack formulation. Section 4 gives some preliminary results. The reduced-load equivalence result is established in Section 5. Section 6 develops the detailed probabilistic arguments involved in deriving the tail asymptotics for the reduced system. In Section 7, we discuss the relationship between the asymptotic regime considered here ("large buffers") and a "many-sources" regime.

**2. Model description.** We first present a detailed model description. We consider a queue with unit capacity (i.e., working at unit speed) fed by several flows indexed by the set $\mathcal{I}$. For any subset $E \subseteq \mathcal{I}$, denote by $A_E(s,t) := \sum_{i \in E} A_i(s,t)$ the aggregate amount of traffic generated by the flows $i \in E$ during the time interval $(s,t]$. Denote by $\rho_E := \sum_{i \in E} \rho_i$ the aggregate traffic intensity of the flows $i \in E$ (as will be specified in detail below). We assume $\rho := \rho_\mathcal{I} < 1$ for stability.

For any $c \geq 0$, $E \subseteq \mathcal{I}$, define $V_E^c(t) := \sup_{0 \leq s \leq t} \{A_E(s,t) - c(t-s)\}$ as the workload at time $t$ in a queue of capacity $c$ fed by the flows $i \in E$ [assuming $V_E^c(0) = 0$]. For $c > \rho_E$, let $\mathbf{V}_E^c$ be a random variable with the limiting distribution of $V_E^c(t)$ for $t \to \infty$. In particular, $V(t) := V_\mathcal{I}^1(t)$ is the total workload, and $\mathbf{V} := \mathbf{V}_\mathcal{I}^1$ is a random variable with the limiting distribution of $V(t)$ for $t \to \infty$.

We assume the flows may be partitioned into two sets: $\mathcal{I}_1$ is the set of "light-tailed" flows; $\mathcal{I}_2$ is the set of "heavy-tailed" flows. For the flows $i \in \mathcal{I}_1$ we make the following assumption.

ASSUMPTION 2.1. For any $c > \rho_{\mathcal{I}_1}$, $\mu > 0$,

$$\lim_{x \to \infty} x^\mu \{ \mathbf{V}_{\mathcal{I}_1}^c > x \} = 0.$$

The above assumption is quite weak; see, for instance, [17] for a very general class of arrival processes satisfying a large-deviations principle (with linear scaling function). However, (superpositions of) On–Off flows of which the activity period has a Weibull distribution satisfy Assumption 2.1 too, as can easily be shown using the bounds in [15] or Section 4.1 of the present paper. Instantaneous renewal input of which the tail of the jump sizes (bursts) is lighter than any power tail is covered by Assumption 2.1 as well.

We assume the flows in $\mathcal{I}_2$ generate traffic according to independent On–Off processes, each alternating between On and Off periods. The Off periods of flow $i$ are generally distributed with mean $1/\lambda_i$. The On periods $\mathbf{A}_i$ have a heavy-tailed distribution $A_i(\cdot)$ with mean $\alpha_i < \infty$. While On, flow $i$ produces traffic at constant rate $r_i$, so the mean burst size is $\alpha_i r_i$. The fraction of time that flow $i$ is On is

$$p_i = \frac{\alpha_i}{1/\lambda_i + \alpha_i} = \frac{\lambda_i \alpha_i}{1 + \lambda_i \alpha_i}.$$

Thus the traffic intensity of flow $i$ is

$$\rho_i := p_i r_i = \frac{\lambda_i \alpha_i r_i}{1 + \lambda_i \alpha_i}.$$

Before stating an important preliminary result, we first introduce some useful notation. For any two real functions $f(\cdot)$ and $g(\cdot)$, we use the notational convention $f(x) \sim g(x)$ to denote $\lim_{x \to \infty} f(x)/g(x) = 1$. Also, we use $f(x) \lesssim g(x)$ to denote $\limsup_{x \to \infty} f(x)/g(x) \leq 1$. Similarly, $f(x) \gtrsim g(x)$ denotes $\liminf_{x \to \infty} f(x)/g(x) \geq 1$. With $X \stackrel{d}{=} Y$ we denote that $X$ and $Y$ have the same distribution.

For any positive stochastic variable $\mathbf{X}$ with distribution function $F(\cdot)$, $\mathbb{E}\{\mathbf{X}\} < \infty$, denote by $F^r(\cdot)$ the distribution function of the residual lifetime of $\mathbf{X}$, that is,

$$F^r(x) := \frac{1}{\mathbb{E}\{\mathbf{X}\}} \int_0^x \mathbb{P}(1 - F(y)) \, dy,$$

and denote by $\mathbf{X}^r$ a stochastic variable with that distribution.

The classes of *long-tailed*, *subexponential*, *regularly varying* and *intermediately regularly varying* distributions are denoted with the symbols $\mathcal{L}$, $\mathcal{S}$, $\mathcal{R}$ and $\mathcal{IR}$, respectively (note that $\mathcal{R} \subset \mathcal{IR} \subset \mathcal{S} \subset \mathcal{L}$). Background on heavy-tailed distributions may be found in [16].

For each flow $i \in \mathcal{I}_2$, we assume that the On period distribution is regularly varying of index $-\nu_i$, that is, $A_i(\cdot) \in \mathcal{R}_{-\nu_i}$ for some $\nu_i > 1$. The next result, which is due to [22], then yields the tail behavior of the workload distribution.

THEOREM 2.1. *If* $A_i^r(\cdot) \in \mathcal{S}$, $\rho_i < c < r_i$, *then*

$$\mathbb{P}\{\mathbf{V}_i^c > x\} \sim (1 - p_i)\frac{\rho_i}{c - \rho_i}\mathbb{P}\Big\{\mathbf{A}_i^r > \frac{x}{r_i - c}\Big\}.$$

**3. Overview of the results.** We now give a broad overview of the main results of the paper. As mentioned in the Introduction, asymptotic bounds in [15] show a sharp dichotomy in the qualitative behavior of $\mathbb{P}\{\mathbf{V} > x\}$, depending on the value of $\rho_{\mathcal{I}_1} + r_{\mathcal{I}_2}$ (i.e., the mean rate of the light-tailed flows plus the peak rate of the heavy-tailed flows) relative to the service rate. In case $\rho_{\mathcal{I}_1} + r_{\mathcal{I}_2} < 1$, the workload has light-tailed characteristics, whereas $\rho_{\mathcal{I}_1} + r_{\mathcal{I}_2} > 1$ implies heavy-tailed characteristics. In the present paper we determine the exact asymptotics of $\mathbb{P}\{\mathbf{V} > x\}$ in the latter case. For the case $\rho_{\mathcal{I}_1} + r_{\mathcal{I}_2} < 1$ (where both the light-tailed and heavy-tailed input determine the workload asymptotics) we refer to [7].

3.1. *Intuitive arguments.* Before formulating our main theorems, we first provide a heuristic derivation of the tail behavior of $\mathbb{P}\{\mathbf{V} > x\}$.

Large-deviations theory suggests that, given that a "rare event" occurs, with overwhelming probability "it happens in the most likely way." In the asymptotic regime considered here ("large buffers"), the most likely way usually consists of a linear build-up of the workload, due to temporary instability of the system. In case of heavy-tailed distributions, the temporary instability typically arises from a "minimal set" of potential causes. The minimal set corresponds to the minimal *number* of causes when these are homogeneous in nature. In general, however, when the potential causes have heterogeneous characteristics, not only the number of them matters, but also their relative likelihood, and their relative contribution to the occurrence of the rare event under consideration.

Translated to our situation, temporary instability is most likely caused by a "minimal set" of flows generating an extreme amount of traffic, while all other flows show roughly average behavior. These considerations give rise to the following characterization of the tail behavior of $\mathbb{P}\{\mathbf{V} > x\}$:

$$\mathbb{P}\{\mathbf{V} > x\} \sim \mathbb{P}\{\mathbf{V}_{S^*}^{c_{S^*}} > x\},$$

with $S^*$ representing the "minimal set," and $c_{S^*} := 1 - \rho_{\mathcal{I}\setminus S^*}$ representing the service rate subtracted by the aggregate traffic intensity of all other flows.

We now introduce some helpful notions in order to formalize the above intuitive arguments. For any subset $S \subseteq \mathcal{I}_2$, define $c_S := 1 - \rho_{\mathcal{I}\setminus S}$ as the service rate subtracted by the aggregate traffic intensity of all other flows $j \in \mathcal{I} \setminus S$. Observe that the stability condition implies $\rho_S < c_S$ for any $S \subseteq \mathcal{I}_2$.

For any subset $S \subseteq \mathcal{I}_2$, denote by $r_S := \sum_{j \in S} r_j$ the aggregate peak rate of the flows $j \in S$. Define $d_S := r_S - c_S = r_S + \rho_{\mathcal{I} \setminus S} - 1$ as the net input rate (i.e., the drift) when all flows in $S$ are On and all other flows show average behavior.

A set $S \subseteq \mathcal{I}_2$ is called (strictly) *critical* if $d_S \geq (>) \, 0$, that is, if

$$r_S + \rho_{\mathcal{I} \setminus S} \geq (>) \, 1.$$

Thus, when all flows in a (strictly) critical set are On, the workload has a (strictly) positive drift. A critical set $S$ is termed *minimally critical* if no proper subset of $S$ is critical, that is, $d_S < \min_{j \in S} \{r_j - \rho_j\}$.

For any subset $S \subseteq \mathcal{I}_2$, denote $\mu_S := \sum_{j \in S} (\nu_j - 1)$. A strictly critical set $S \subseteq \mathcal{I}_2$ is said to be (weakly) *dominant* if $\mu_S < (\leq) \, \mu_U$ for any other critical set $U \subseteq \mathcal{I}_2$. Observe that for a set $S \subseteq \mathcal{I}_2$ to be dominant, it must be minimally critical (because otherwise the defining property would be violated for any critical subset $U \subset S$).

The quantity $\mu_S$ may be interpreted as a measure for the "cost" associated with a temporary drift $d_S$: the probability of all flows in $S$ being On for a time of the order $x$ in steady state is roughly equal to $x^{-\mu_S}$. Thus, a set $S$ is (weakly) dominant if the flows in $S$ being On causes the drift to be positive in the cheapest possible way.

In case of light-tailed distributions, the cost minimization is usually not so simple; one then also needs to consider how long a certain positive drift must be maintained in order for a given workload level $x$ to be reached. This issue does not arise in case of regularly varying On periods, since $\mathbb{P}\{\mathbf{A}_i^r > ax\}$ is of the same order of magnitude (up to a constant) as $\mathbb{P}\{\mathbf{A}_i^r > x\}$ for any constant $a > 1$. This implies that the value of the temporary drift is not relevant as long as it is positive.

Note that these heuristic arguments clearly do not hold for other subexponential distributions, such as the lognormal and Weibull distribution. In this case, one has $\mathbb{P}\{\mathbf{A}_i^r > ax\} = o(\mathbb{P}\{\mathbf{A}_i^r > x\})$, if $a > 1$.

3.2. *Tail behavior of the workload distribution.* We now state our main theorem.

THEOREM 3.1 (Reduced-load equivalence). *Suppose the set of flows $S^* \subseteq \mathcal{I}_2$ is dominant. If $A_j(\cdot) \in \mathcal{R}$ for all $j \in \mathcal{I}_2$, then*

$$(3.1) \qquad \mathbb{P}\{\mathbf{V} > x\} \sim \mathbb{P}\{\mathbf{V}_{S^*}^{c_{S^*}} > x\},$$

*with*

$$(3.2) \qquad \mathbb{P}\{\mathbf{V}_{S^*}^{c_{S^*}} > x\} \sim \left( \prod_{j \in S^*} p_j \right) \sum_{\mathcal{J}_0 \subseteq S^*} P_{\mathcal{J}_0}(x),$$

*where $P_{\mathcal{J}_0}(x)$ is given by (with $\mathcal{J}_1 = S^* \setminus \mathcal{J}_0$, and $d_{S^*} = r_{S^*} - c_{S^*}$ as defined earlier)*

$$P_{\mathcal{J}_0}(x) = \frac{1}{\prod_{i \in \mathcal{J}_1} \mathbb{E}\{\mathbf{A}_i\}}$$

(3.3)
$$\times \int_{y_i \in (0,\infty), i \in \mathcal{J}_1} \prod_{i \in \mathcal{J}_1} \mathbb{P}\left\{ d_{S^*} \mathbf{A}_i > \sum_{j \in \mathcal{J}_1} y_j(r_j - \rho_j) - d_{S^*} y_i + x \right\}$$

$$\times \prod_{i \in \mathcal{J}_0} \mathbb{P}\left\{ d_{S^*} \mathbf{A}_i^r > \sum_{j \in \mathcal{J}_1} y_j(r_j - \rho_j) + x \right\} \prod_{i \in \mathcal{J}_1} dy_i.$$

*In particular, $\mathbb{P}\{\mathbf{V} > x\}$ and $P_{\mathcal{J}_0}(x)$ are regularly varying of index $-\mu_{S^*} = -\sum_{j \in S^*}(\nu_j - 1)$.*

The proof of the above theorem may be found in Section 5.1 [(3.1)] and Section 6 [(3.2) and (3.3) and the regular variation property].

Note that in case the reduced system consists of just a single flow, that is, $S^* = \{i^*\}$, the tail asymptotics follow directly from Theorem 2.1. This is in fact the reduced-load equivalence established in [1] (under somewhat weaker distributional assumptions). Note that in this case the right-hand side of (3.2) takes the form $p_{i^*}[P_{\varnothing}(x) + P_{i^*}(x)]$, with

$$P_{i^*}(x) = \mathbb{P}\left\{ \mathbf{A}_{i^*}^r > \frac{x}{r_{i^*} - c_{i^*}} \right\},$$

and (after a straightforward calculation)

$$P_{\varnothing}(x) = \frac{r_{i^*} - c_{i^*}}{c_{i^*} - \rho_{i^*}} \mathbb{P}\left\{ \mathbf{A}_{i^*}^r > \frac{x}{r_{i^*} - c_{i^*}} \right\},$$

so that

$$p_{i^*}[P_{\varnothing}(x) + P_{i^*}(x)] = (1 - p_{i^*}) \frac{\rho_{i^*}}{c_{i^*} - \rho_{i^*}} \mathbb{P}\left\{ \mathbf{A}_{i^*}^r > \frac{x}{r_{i^*} - c_{i^*}} \right\},$$

which is consistent with Theorem 2.1.

In case the reduced system consists of several flows, the tail asymptotics cannot be obtained from known results. In fact, the analysis of the reduced system then poses a major challenge because of the rather subtle mechanics involved in reaching a large workload level. By definition, though, the reduced system has the special feature that all flows must be On for the drift in the workload to be positive, that is, $r_{S^*} - \min_{j \in S^*}\{r_j - \rho_j\} < c_{S^*} < r_{S^*}$. In Section 6 we determine the exact asymptotics for systems satisfying this property, yielding the integral expression given in Theorem 3.1.

3.3. *Knapsack formulation for determining a dominant set.* We now describe how a dominant set may be determined from a simple knapsack formulation. Recall that the On period distributions of the flows $i \in \mathcal{I}_2$ are regularly varying of index $-\nu_i$.

For a strictly critical set $S \subseteq \mathcal{I}_2$ to be dominant, it must necessarily solve the optimization problem

$$\min_{S \subseteq \mathcal{I}_2} \sum_{j \in S} (\nu_j - 1)$$

$$\text{sub} \quad \sum_{j \in S} r_j + \sum_{j \in \mathcal{I}_2 \setminus S} \rho_j > 1 - \rho_{\mathcal{I}_1}.$$

Note that the constraint is equivalent to $d_S > 0$. If we define $\theta_i := r_i - \rho_i$ for all $i \in \mathcal{I}_2$, then the above problem may be expressed in the standard knapsack form as

$$\max_{U \subseteq \mathcal{I}_2} \sum_{j \in U} (\nu_j - 1)$$

$$\text{sub} \quad \sum_{j \in U} \theta_j \leq \rho_{\mathcal{I}_1} + r_{\mathcal{I}_2} - 1 - \varepsilon,$$

with $U = \mathcal{I}_2 \setminus S$ and $\varepsilon$ some small positive number. The above problem may not always have a unique solution. In case it does, the corresponding set $S$ is dominant, except for the case when some set $T$ exists which is critical but not strictly critical (i.e., $r_T + \rho_{I \setminus T} = 1$), with $\mu_T \leq \mu_S$ (see the definition of a dominant set). Although intriguing, this "critical case" is not further considered in the present paper. In this case, the temporary drift may be *zero* for a long period of time during the path to overflow.

In case the knapsack problem has several solutions, the corresponding sets are weakly dominant (except for the critical case again). The next theorem extends the reduced-load equivalence to the case of weakly dominant sets.

THEOREM 3.2 (Generalized reduced-load equivalence; weakly dominant sets). *Let $\Upsilon \subseteq 2^{\mathcal{I}_2}$ be the collection of all weakly dominant sets. If $A_j(\cdot) \in \mathcal{R}$ for all $j \in S$, $S \in \Upsilon$, then*

$$(3.4) \qquad \mathbb{P}\{\mathbf{V} > x\} \sim \sum_{S \in \Upsilon} \mathbb{P}\{\mathbf{V}_S^{c_S} > x\},$$

*with $\mathbb{P}\{\mathbf{V}_S^{c_S} > x\}$ as in* (3.2), (3.3).

3.4. *Homogeneous On–Off flows.* We briefly consider the case of homogeneous On–Off flows as an important special case with weakly dominant sets. Assume that the flows $i \in \mathcal{I}_2$ have identical characteristics. With some minor abuse of notation, let $A(\cdot) := A_i(\cdot)$, $\nu := \nu_i$, $\rho := \rho_i$, $r := r_i$, $p_i \equiv p$.

Define $N^* := \arg\min\{N : Nr + (|\mathcal{I}_2| - N)\rho > 1 - \rho_{\mathcal{I}_1}\}$. (Observe that the assumption $\rho_{\mathcal{I}_1} + r_{\mathcal{I}_2} > 1$ ensures $N^* \le |\mathcal{I}_2|$.) To exclude the critical case, assume that $(N^* - 1)r + (|\mathcal{I}_2| - N^* + 1)\rho < 1 - \rho_{\mathcal{I}_1}$, so that the drift remains negative (and cannot be zero) when only $N^* - 1$ flows are On.

COROLLARY 3.1. *If $A(\cdot) \in \mathcal{R}$, then*

$$\mathbb{P}\{\mathbf{V} > x\} \sim \binom{|\mathcal{I}_2|}{N^*} \mathbb{P}\{\bar{\mathbf{V}} > x\},$$

*with*

$$\mathbb{P}\{\bar{\mathbf{V}} > x\} \sim p^{N^*} \sum_{n=0}^{N^*} \binom{N^*}{n} P_{\{1,\dots,n\}}(x),$$

*where $P_{\{1,\dots,n\}}(x)$ is given by (3.3). In particular, $\mathbb{P}\{\mathbf{V} > x\}$ and $P_{\{1,\dots,n\}}(x)$ are regularly varying of index $-N^*(\nu - 1)$.*

3.5. *$K$ heterogeneous classes.* We finally consider the important special case where each On–Off flow in $\mathcal{I}_2$ belongs to one of $K$ heterogeneous classes. We will show how an approximate solution to the knapsack problem may be obtained using a simple index rule. The approximation is in fact asymptotically exact in the many-sources regime.

Specifically, consider the superposition of $n$ On–Off flows, each belonging to one of $K$ heterogeneous classes. Let $a_k$ be the fraction of flows of class $k \in \{1, \dots, K\}$, with peak rate $r_k$, mean rate $\rho_k$, and an On period distribution which is regularly varying of index $-\nu_k$. Let the service rate be $n$ (instead of 1), and let $\mathbf{V}^{(n)}$ be the stationary workload. The knapsack problem then takes the form

$$\min_{n_k \in \{0,\dots,na_k\}} \sum_{k=1}^{K} n_k(\nu_k - 1)$$

$$\text{sub} \quad \sum_{k=1}^{K} n_k r_k + \sum_{k=1}^{K} (na_k - n_k)\rho_k > n.$$

Unfortunately, the above problem cannot be easily solved due to the integrality constraints. Intuitively, however, one may expect that as $n$ grows large, the integrality constraints should have a negligible effect, so that a continuous relaxation with $n_k \in [0, na_k]$ should give a good approximate solution.

This relaxation may be solved using a simple index rule. Index the $K$ classes in nondecreasing order of the ratios

$$\gamma_k := (\nu_k - 1)/(r_k - \rho_k).$$

For any $k \in \{1, \ldots, K\}$, define $\sigma_k := \sum_{m=1}^{k-1} a_m r_m + \sum_{m=k}^{K} a_m \rho_m$. Determine the (unique) index $l$ such that $1 \in (\sigma_{l-1}, \sigma_l]$. Then take $n_k^* = na_k$ for all classes $k < l$, $n_k^* = 0$ for all classes $k > l$, and $n_l^* = n(1 - \sigma_{l-1})/(r_l - \rho_l)$.

This yields the (crude) approximation

$$(3.5) \qquad \mathbb{P}\{\mathbf{V}^{(n)} > x\} \approx x^{-n\mu},$$

with $\mu := \sum_{k=1}^{l-1} a_k(\nu_k - 1) + (1 - \sigma_{l-1})\gamma_l$. In Section 7 we prove that the above approximation is logarithmically exact in the many-sources regime. In particular, one may show that the limits for $x \to \infty$ and $n \to \infty$ commute if one considers logarithmic asymptotics.

THEOREM 3.3 (Robustness of logarithmic asymptotics).

$$\lim_{n \to \infty} \lim_{x \to \infty} \frac{1}{n} \frac{\log \mathbb{P}\{\mathbf{V}^{(n)} > nx\}}{\log x} = \lim_{x \to \infty} \lim_{n \to \infty} \frac{1}{n} \frac{\log \mathbb{P}\{\mathbf{V}^{(n)} > nx\}}{\log x}.$$

The proof of the above theorem may be found in Section 7. Although logarithmically exact, the approximation (3.5) may not be appropriate from a practical perspective. In particular, it is shown in Section 7 that an analogue of Theorem 3.3 cannot hold if one considers exact asymptotics. This "negative" result is reminiscent of a phenomenon occurring in heavy-traffic theory where two limiting regimes lead either to stable Lévy motion or to fractional Brownian motion; see, for example, [31] and references therein.

**4. Preliminary results.** In this section we collect some preliminary results which will be used in later sections.

4.1. *Bounds.* We first derive some simple bounds for the workload distribution $\mathbb{P}\{\mathbf{V}_S^c > x\}$ for subsets $S \subseteq \mathcal{I}_2$.

For any subset $S \subseteq \mathcal{I}_2$, $c < r_S$, define

$$P_S^c(x) := \prod_{j \in S} p_j \mathbb{P}\left\{\mathbf{A}_j^r > \frac{x}{r_S - c}\right\}.$$

The next lemma gives a lower bound for $\mathbb{P}\{\mathbf{V}_S^c > x\}$ which may also be found in [12].

LEMMA 4.1. *Let $S \subseteq \mathcal{I}_2$. For $c < r_S$,*

$$\mathbb{P}\{\mathbf{V}_S^c > x\} \geq P_S^c(x).$$

PROOF. Consider the event that at some arbitrary time $t$ all flows $j \in S$ have been On since time $t - \frac{x}{r_S - c}$ or longer. This event occurs with probability $P_s^c(x)$, and implies that the workload at time $t$ is larger than $\frac{r_S x}{r_S - c} - \frac{cx}{r_S - c} = x$. $\quad\square$

For any subset $S \subseteq \mathcal{I}_2$, $c < r_S$, define

$$K_S^c := \prod_{j \in S} \frac{r_j - \rho_j}{r_j - \rho_j + c - r_S}.$$

The next lemma establishes an asymptotic upper bound for $\mathbb{P}\{\mathbf{V}_S^c > x\}$ for the case where $S$ is a minimally critical set with respect to the capacity $c$.

LEMMA 4.2. *Let* $S \subseteq \mathcal{I}_2$. *If* $c \in (r_S - \min_{j \in S}\{r_j - \rho_j\}, r_S)$, *and* $A_j^r(\cdot) \in \mathcal{S}$ *for all* $j \in S$, *then*

$$\mathbb{P}\{\mathbf{V}_S^c > x\} \lesssim K_S^c P_S^c(x).$$

PROOF. For any $i \in S$, denote $d_i := c - r_S + r_i$. Observe that $d_i > \rho_i$ since $c > r_S - (r_i - \rho_i)$. We apply the usual technique to obtain an upper bound: split the capacity. Formally, we have the sample-path upper bound

$$(4.1) \qquad V_S^c(t) \leq V_i^{d_i}(t) + V_{S\setminus\{i\}}^{r_{S\setminus\{i\}}}(t) = V_i^{d_i}(t)$$

for all $i \in S$.

In the stationary regime, using Theorem 2.1,

$$\begin{aligned}
\mathbb{P}\{\mathbf{V}_S^c > x\} &\leq \mathbb{P}\{\mathbf{V}_j^{d_j} > x \text{ for all } j \in S\} \\
&= \prod_{j \in S} \mathbb{P}\{\mathbf{V}_j^{d_j} > x\} \\
&\sim \prod_{j \in S}(1 - p_j)\frac{\rho_j}{d_j - \rho_j}\mathbb{P}\left\{\mathbf{A}_j^r > \frac{x}{r_j - d_j}\right\} \\
&= \prod_{j \in S} p_j \frac{r_j - \rho_j}{r_j - \rho_j + c - r_S}\mathbb{P}\left\{\mathbf{A}_j^r > \frac{x}{r_S - c}\right\} \\
&= K_S^c P_S^c(x). \qquad\qquad\qquad\qquad\qquad\quad \square
\end{aligned}$$

COROLLARY 4.1. *Let* $S \subseteq \mathcal{I}_2$. *If* $c \in (r_S - \min_{j \in S}\{r_j - \rho_j\}, r_S)$, *and* $A_j^r(\cdot) \in \mathcal{S}$ *for all* $j \in S$, *then*

$$P_S^c(x) \leq \mathbb{P}\{\mathbf{V}_S^c > x\} \lesssim K_S^c P_S^c(x).$$

PROOF. The proof follows directly by combining Lemmas 4.1 and 4.2. $\square$

COROLLARY 4.2. *Let $S \subseteq \mathcal{I}_2$. If $A_j^r(\cdot) \in \mathcal{IR}$ for all $j \in S$, then, for any closed interval $T \subseteq (r_S - \min_{j \in S}\{r_j - \rho_j\}, r_S)$, there exist constants $K^{(1)}$, $K^{(2)}$ independent of $c$, such that for all $c \in T$,*

$$K^{(1)} P_S(x) \lesssim \mathbb{P}\{\mathbf{V}_S^c > x\} \lesssim K^{(2)} P_S(x),$$

*with*

$$P_S(x) := \prod_{j \in S} \mathbb{P}\{\mathbf{A}_j^r > x\}.$$

PROOF. The statement follows directly from Corollary 4.1 and the fact that $A_j^r(\cdot) \in \mathcal{IR} \subset \mathcal{S}$ for all $j \in S$ when observing that $A_j^r(\cdot) \in \mathcal{IR}$, $j \in S$ implies that

$$\limsup_{x \to \infty} \frac{P_S^{c_1}(x)}{P_S^{c_2}(x)} < \infty,$$

if $c_1, c_2 \in T$. □

We now derive some general bounds for the total workload distribution $\mathbb{P}\{\mathbf{V} > x\}$ which will be crucial in establishing the reduced-load equivalence.

For any $c \geq 0$, $E \subseteq \mathcal{I}$, define $Z_E^c(t) := \sup_{0 \leq s \leq t}\{c(t-s) - A_E(s,t)\}$. For $c < \rho_E$, let $\mathbf{Z}_E^c$ be a random variable with the limiting distribution of $Z_E^c(t)$ for $t \to \infty$. Let $\Omega \subseteq 2^{\mathcal{I}_2}$ be the collection of all minimally critical sets.

We first present a lower bound. The idea is as follows: $\mathbf{V}_E^{c_E}$ being large for some minimally critical set $E \in \Lambda$ basically implies that $\mathbf{V}$ must be large too, unless the other flows $j \notin E$ persist in below-average behavior. Excluding such below-average behavior (reflected in large values of $\mathbf{Z}_{\mathcal{I}\backslash E}^c$) from the event $\{\mathbf{V} > x\}$ yields the following lower bound for $\mathbb{P}\{\mathbf{V} > x\}$.

LEMMA 4.3. *Let $\Lambda \subseteq \Omega$. Then for any $\delta > 0$ and $y \geq 0$,*

$$\mathbb{P}\{\mathbf{V} > x\} \geq \sum_{E \in \Lambda} \mathbb{P}\{\mathbf{V}_E^{c_E+\delta} > x+y\} \mathbb{P}\{\mathbf{Z}_{\mathcal{I}\backslash E}^{\rho_{\mathcal{I}\backslash E}-\delta} \leq y\}$$

$$- \sum_{E_1, E_2 \in \Lambda, E_1 \neq E_2} \prod_{j \in E_1 \cup E_2} \mathbb{P}\left\{\mathbf{V}_j^{\rho_j+\delta} > x\right\}.$$

PROOF. Sample-path wise,

$$V(t) = \sup_{0 \leq s \leq t}\{A(s,t) - (t-s)\}$$

$$= \sup_{0 \leq s \leq t}\{A_E(s,t) + A_{\mathcal{I}\backslash E}(s,t) - (c_E+\delta)(t-s)$$

$$- (\rho_{\mathcal{I}\backslash E} - \delta)(t-s)\}$$

$$\geq \sup_{0 \leq s \leq t} \{A_E(s,t) - (c_E + \delta)(t-s)\}$$

$$+ \inf_{0 \leq s \leq t} \{A_{\mathcal{I} \setminus E}(s,t) - (\rho_{\mathcal{I} \setminus E} - \delta)(t-s)\}$$

$$= \sup_{0 \leq s \leq t} \{A_E(s,t) - (c_E + \delta)(t-s)\}$$

$$- \sup_{0 \leq s \leq t} \{(\rho_{\mathcal{I} \setminus E} - \delta)(t-s) - A_{\mathcal{I} \setminus E}(s,t)\}$$

$$= V_E^{c_E + \delta}(t) - Z_{\mathcal{I} \setminus E}^{\rho_{\mathcal{I} \setminus E} - \delta}(t)$$

for all $E \in \Lambda$.

In the stationary regime, for any $\delta > 0$ and $y \geq 0$, using the independence of $\mathbf{V}_E^{c_E + \delta}$ and $\mathbf{Z}_{\mathcal{I} \setminus E}^{\rho_{\mathcal{I} \setminus E} - \delta}$,

$$\mathbb{P}\{\mathbf{V} > x\}$$

$$\geq \mathbb{P}\{\mathbf{V}_E^{c_E + \delta} - \mathbf{Z}_{\mathcal{I} \setminus E}^{\rho_{\mathcal{I} \setminus E} - \delta} > x \text{ for some } E \in \Lambda\}$$

$$\geq \mathbb{P}\{\mathbf{V}_E^{c_E + \delta} > x + y, \mathbf{Z}_{\mathcal{I} \setminus E}^{\rho_{\mathcal{I} \setminus E} - \delta} \leq y \text{ for some } E \in \Lambda\}$$

(4.2)
$$\geq \mathbb{P}\{\mathbf{V}_E^{c_E + \delta} > x + y, \mathbf{Z}_{\mathcal{I} \setminus E}^{\rho_{\mathcal{I} \setminus E} - \delta} \leq y \text{ for exactly one } E \in \Lambda\}$$

$$= \sum_{E \in \Lambda} \mathbb{P}\{\mathbf{V}_E^{c_E + \delta} > x + y, \mathbf{Z}_{\mathcal{I} \setminus E}^{\rho_{\mathcal{I} \setminus E} - \delta} \leq y\}$$

$$- \sum_{E_1, E_2 \in \Lambda, E_1 \neq E_2} \mathbb{P}\{\mathbf{V}_{E_1}^{c_{E_1} + \delta} > x + y, \mathbf{Z}_{\mathcal{I} \setminus E_1}^{\rho_{\mathcal{I} \setminus E_1} - \delta} \leq y,$$

$$\mathbf{V}_{E_2}^{c_{E_2} + \delta} > x + y, \mathbf{Z}_{\mathcal{I} \setminus E_2}^{\rho_{\mathcal{I} \setminus E_2} - \delta} \leq y\}$$

$$\geq \sum_{E \in \Lambda} \mathbb{P}\{\mathbf{V}_E^{c_E + \delta} > x + y\}\mathbb{P}\{\mathbf{Z}_{\mathcal{I} \setminus E}^{\rho_{\mathcal{I} \setminus E} - \delta} \leq y\}$$

$$- \sum_{E_1, E_2 \in \Lambda, E_1 \neq E_2} \mathbb{P}\{\mathbf{V}_{E_1}^{c_{E_1} + \delta} > x, \mathbf{V}_{E_2}^{c_{E_2} + \delta} > x\}.$$

As in (4.1),

(4.3)
$$V_E^{c_E + \delta}(t) \leq V_i^{c_E - r_{E \setminus \{i\}} + \delta}(t) + V_{E \setminus \{i\}}^{r_{E \setminus \{i\}}}(t)$$

$$= V_i^{c_E - r_{E \setminus \{i\}} + \delta}(t)$$

for all $i \in E$.

Note that $c_E - r_{E \setminus \{i\}} > \rho_i$ for all $i \in E$, $E \in \Lambda$, since $E$ is minimally critical.

Hence,

$$V_E^{c_E+\delta}(t) \leq V_i^{\rho_i+\delta}(t)$$

for all $i \in E$, $E \in \Lambda$.

Thus,

$$\mathbb{P}\{\mathbf{V}_{E_1}^{c_{E_1}+\delta} > x, \mathbf{V}_{E_2}^{c_{E_2}+\delta} > x\}$$

(4.4)
$$\leq \mathbb{P}\{\mathbf{V}_j^{\rho_j+\delta} > x \text{ for all } j \in E_1, \mathbf{V}_j^{\rho_j+\delta} > x \text{ for all } j \in E_2\}$$

$$= \mathbb{P}\{\mathbf{V}_j^{\rho_j+\delta} > x \text{ for all } j \in E_1 \cup E_2\}$$

$$= \prod_{j \in E_1 \cup E_2} \mathbb{P}\{\mathbf{V}_j^{\rho_j+\delta} > x\}.$$

Substituting (4.4) into (4.2) completes the proof. $\square$

We now provide a corresponding upper bound, which is somewhat more involved. The idea is as follows: $\mathbf{V}$ being large essentially means that $\mathbf{V}_E^{c_E}$ must be large for some minimally critical set $E \in \Lambda$ too, unless the other flows $j \notin E$ exhibit above-average behavior. Extending the event $\{\mathbf{V} > x\}$ with possible above-average behavior of the flows $j \notin E$ (manifesting itself in large values of $\mathbf{V}_{\mathcal{I} \setminus E}^{\rho_{\mathcal{I} \setminus E}+\delta}$) leads to the following upper bound for $\mathbb{P}\{\mathbf{V} > x\}$.

LEMMA 4.4. *Let $\Lambda \subseteq \Omega$. Then for any $\delta, \varepsilon > 0$ sufficiently small and $y$,*

$$\mathbb{P}\{\mathbf{V} > x\} \leq \sum_{E \in \Lambda} \mathbb{P}\{\mathbf{V}_E^{c_E-\delta} > x-y\} + \mathbb{P}\{\mathbf{V}_{\mathcal{I}_1}^{\rho_{\mathcal{I}_1}+\varepsilon} > x/\mathcal{N}\}$$

$$+ \sum_{E \in \Lambda} \mathbb{P}\{\mathbf{V}_{\mathcal{I} \setminus E}^{\rho_{\mathcal{I} \setminus E}+\delta} > y\} \prod_{j \in E} \mathbb{P}\{\mathbf{V}_j^{\rho_j+\varepsilon} > x/\mathcal{N}\}$$

$$+ \sum_{E \in \Omega \setminus \Lambda} \prod_{j \in E} \mathbb{P}\{\mathbf{V}_j^{\rho_j+\varepsilon} > x/\mathcal{N}\},$$

*with $\mathcal{N} := |\mathcal{I}|$ denoting the total number of flows.*

PROOF. As before, we divide the capacity to obtain the sample-path upper bound

$$V(t) \leq V_E^{c_E-\delta}(t) + V_{\mathcal{I} \setminus E}^{\rho_{\mathcal{I} \setminus E}+\delta}(t)$$

for all $E \in \Lambda$.

In addition, for $\varepsilon > 0$ sufficiently small, $V(t) > x$ implies $V_{\mathcal{I}_1}^{\rho_{\mathcal{I}_1}+\varepsilon}(t) > x/\mathcal{N}$, or there exists a minimally critical set $S \in \Omega$ such that $V_j^{\rho_j+\varepsilon}(t) > x/\mathcal{N}$ for all $j \in S$.

This may be seen as follows: suppose that it were not the case, that is, $V_{\mathcal{I}_1}^{\rho_{\mathcal{I}_1}+\varepsilon}(t) \leq x/\mathcal{N}$, and for every minimally critical set $S \in \Omega$ there exists a $j$ (depending on $S$) such that $V_j^{\rho_j+\varepsilon}(t) \leq x/\mathcal{N}$. Then the set $\mathcal{J}(t) := \{j \in \mathcal{I}_2 : V_j^{\rho_j+\varepsilon}(t) > x/\mathcal{N}\}$ does not contain any minimally critical set, hence $r_{\mathcal{J}(t)} + \rho_{\mathcal{I}\setminus\mathcal{J}(t)} < 1$. This means that $\rho_{\mathcal{I}\setminus\mathcal{J}(t)} + \mathcal{N}\varepsilon \leq 1 - r_{\mathcal{J}(t)}$ for $\varepsilon > 0$ sufficiently small. Thus, noting that $\rho_{\mathcal{I}\setminus\mathcal{J}(t)} = \rho_{\mathcal{I}_1} + \rho_{\mathcal{I}_2\setminus\mathcal{J}(t)}$,

$$V(t) \leq V_{\mathcal{J}(t)}^{r_{\mathcal{J}(t)}}(t) + V_{\mathcal{I}\setminus\mathcal{J}(t)}^{1-r_{\mathcal{J}(t)}}(t)$$

$$= V_{\mathcal{I}\setminus\mathcal{J}(t)}^{1-r_{\mathcal{J}(t)}}(t)$$

$$\leq V_{\mathcal{I}\setminus\mathcal{J}(t)}^{\rho_{\mathcal{I}\setminus\mathcal{J}(t)}+\mathcal{N}\varepsilon}(t)$$

$$\leq V_{\mathcal{I}_1}^{\rho_{\mathcal{I}_1}+\varepsilon}(t) + \sum_{j\in\mathcal{I}_2\setminus\mathcal{J}(t)} V_j^{\rho_j+\varepsilon}(t)$$

$$\leq |\mathcal{I}\setminus\mathcal{J}(t)|\, x/\mathcal{N}$$

$$\leq x,$$

contradicting the initial supposition.

In the stationary regime, for any $\delta, \varepsilon > 0$ sufficiently small and $y$, using independence

$$\mathbb{P}\{\mathbf{V} > x\}$$

$$\leq \mathbb{P}\{\mathbf{V}_E^{c_E-\delta} + \mathbf{V}_{\mathcal{I}\setminus E}^{\rho_{\mathcal{I}\setminus E}+\delta} > x \text{ for all } E \in \Lambda,$$

$$\mathbf{V}_{\mathcal{I}_1}^{\rho_{\mathcal{I}_1}+\varepsilon} > x/\mathcal{N} \text{ or } \mathbf{V}_j^{\rho_j+\varepsilon} > x/\mathcal{N} \text{ for all } j \in S \text{ for some } S \in \Omega\}$$

$$\leq \mathbb{P}\{\mathbf{V}_E^{c_E-\delta} > x - y \text{ or } \mathbf{V}_{\mathcal{I}\setminus E}^{\rho_{\mathcal{I}\setminus E}+\delta} > y \text{ for all } E \in \Lambda,$$

$$\mathbf{V}_{\mathcal{I}_1}^{\rho_{\mathcal{I}_1}+\varepsilon} > x/\mathcal{N} \text{ or } \mathbf{V}_j^{\rho_j+\varepsilon} > x/\mathcal{N} \text{ for all } j \in S \text{ for some } S \in \Omega\}$$

$$\leq \sum_{E\in\Lambda} \mathbb{P}\{\mathbf{V}_E^{c_E-\delta} > x - y\} + \mathbb{P}\{\mathbf{V}_{\mathcal{I}_1}^{\rho_{\mathcal{I}_1}+\varepsilon} > x/\mathcal{N}\}$$

$$+ \sum_{S\in\Omega} \mathbb{P}\{\mathbf{V}_j^{\rho_j+\varepsilon} > x/\mathcal{N} \text{ for all } j \in S, \mathbf{V}_{\mathcal{I}\setminus E}^{\rho_{\mathcal{I}\setminus E}+\delta} > y \text{ for all } E \in \Lambda\}$$

$$\leq \sum_{E\in\Lambda} \mathbb{P}\{\mathbf{V}_E^{c_E-\delta} > x - y\} + \mathbb{P}\{\mathbf{V}_j^{\rho_j+\varepsilon} > x/\mathcal{N}\}$$

$$+ \sum_{E\in\Lambda} \mathbb{P}\{\mathbf{V}_j^{\rho_j+\varepsilon} > x/\mathcal{N} \text{ for all } j \in E, \mathbf{V}_{\mathcal{I}\setminus E}^{\rho_{\mathcal{I}\setminus E}+\delta} > y\}$$

$$+ \sum_{E\in\Omega\setminus\Lambda} \mathbb{P}\{\mathbf{V}_j^{\rho_j+\varepsilon} > x/\mathcal{N} \text{ for all } j \in E\}$$

$$\leq \sum_{E \in \Lambda} \mathbb{P}\{\mathbf{V}_E^{c_E - \delta} > x - y\} + \mathbb{P}\{\mathbf{V}_j^{\rho_j + \varepsilon} > x/\mathcal{N}\}$$

$$+ \sum_{E \in \Lambda} \mathbb{P}\{\mathbf{V}_{\mathcal{I}\backslash E}^{\rho_{\mathcal{I}\backslash E} + \delta} > y\} \prod_{j \in E} \mathbb{P}\{\mathbf{V}_j^{\rho_j + \varepsilon} > x/\mathcal{N}\}$$

$$+ \sum_{E \in \Omega \backslash \Lambda} \prod_{j \in E} \mathbb{P}\{\mathbf{V}_j^{\rho_j + \varepsilon} > x/\mathcal{N}\}. \qquad \Box$$

4.2. *Stationary workload representation.* We now give a convenient representation for the stationary workload $\mathbf{V}_E^c$, with $E \subseteq \mathcal{I}_2$ an arbitrary set of heavy-tailed On–Off flows. We start from the definition $V_E^c(t) := \sup_{0 \leq s \leq t}\{A_E(s, t) - c(t - s)\}$ [assuming $V_E^c(0) = 0$]. Since the process $A_E(\cdot, \cdot)$ has stationary and reversible increments, we have

$$\sup_{0 \leq s \leq t}\{A_E(s, t) - c(t - s)\} \stackrel{d}{=} \sup_{0 \leq s \leq t}\{A_E(0, s) - cs\}.$$

In the sequel, we simply use the latter expression as the *definition* of $V_E^c(t)$. Accordingly, for $c > \rho_E$, the stationary workload as $t \to \infty$ may be represented as

$$\mathbf{V}_E^c := \sup_{t \geq 0}\{A_E(0, t) - ct\}.$$

Explicit constructions of $A_i(0, t)$ (satisfying the stationarity condition) may be found in [15, 19]. For completeness, we review the construction in [19], which will be extensively used in Section 6.

Let $\{\mathbf{A}_{im}, m \geq 0\}$ be a sequence of i.i.d. random variables representing On periods of flow $i$. Similarly, let $\{\mathbf{U}_{im}, m \geq 1\}$ be Off periods. Define three additional random variables $\mathbf{A}_{i0}^r$, $\mathbf{U}_{i0}^r$ and $\mathbf{B}_i$ such that $\mathbf{A}_{i0}^r \stackrel{d}{=} \mathbf{A}_i^r$, $\mathbf{U}_{i0}^r \stackrel{d}{=} \mathbf{U}_i^r$ and

$$\mathbb{P}\{\mathbf{B}_i = 1\} = \frac{\mathbb{E}\{\mathbf{A}_{i1}\}}{\mathbb{E}\{\mathbf{A}_{i1}\} + \mathbb{E}\{\mathbf{U}_{i1}\}} = 1 - \mathbb{P}\{\mathbf{B}_i = 0\}.$$

Note that $\mathbf{B}_i = 1$ corresponds to flow $i$ being On (in stationarity).

To obtain a stationary alternating renewal process, we define the delay random variable $\mathbf{D}_{i0}$ by

$$\mathbf{D}_{i0} = \mathbf{B}_i \mathbf{A}_{i0}^r + (1 - \mathbf{B}_i)(\mathbf{U}_{i0}^r + \mathbf{A}_{i0}).$$

Then the delayed renewal sequence

$$\{\mathbf{Z}_{in}, n \geq 0\} = \left\{\mathbf{D}_{i0}, \mathbf{D}_{i0} + \sum_{m=1}^{n}(\mathbf{U}_{im} + \mathbf{A}_{im}), n \geq 1\right\}$$

is stationary.

Next, we define the process $\{J_i(t), t \geq 0\}$ as follows. $J_i(t)$ is the indicator of the event that flow $i$ is On at time $t$. Formally, we have

$$J_i(t) = \mathbf{B}_i \mathbb{1}_{\{t < \mathbf{A}_{i0}^r\}} + (1 - \mathbf{B}_i) \mathbb{1}_{\{\mathbf{U}_{i0}^r \leq t < \mathbf{U}_{i0}^r + \mathbf{A}_{i0}\}}$$

$$+ \sum_{n=0}^{\infty} \mathbb{1}_{\{\mathbf{Z}_{in} + \mathbf{U}_{i,n+1} \leq t < \mathbf{Z}_{i,n+1}\}}.$$

The On–Off process $\{J_i(t), t \geq 0\}$ is strictly stationary; see [19], Theorem 2.1. The process $\{A_i(0, t), t \geq 0\}$ is defined by

$$A_i(0, t) := r_i \int_0^t J_i(u) \, du.$$

Finally, note that the number of elapsed Off periods during $[0, t]$ which started after time 0 is given by

$$(4.5) \qquad N_i^A(t) := \max\{n : \mathbf{Z}_{i,n-1} + \mathbf{U}_{in} \leq t\}.$$

We conclude this section with the following useful lemma.

LEMMA 4.5.  *Let $S \subseteq \mathcal{I}_2$. If $A_j(\cdot) \in \mathcal{R}$ for all $j \in S$ and $c \in (r_S - \min_{j \in S}\{r_j - \rho_j\}, r_S)$, then*

$$\lim_{M \to \infty} \limsup_{x \to \infty} \frac{\mathbb{P}\{\sup_{t \geq Mx}\{A_S(0, t) - (c - \varepsilon)t\} > x\}}{\mathbb{P}\{\mathbf{V}_S^c > x\}} = 0,$$

*for any $\varepsilon \in [0, r_S - c)$.*

PROOF.  For $t \geq Mx$, write

$$A_S(0, t) - (c - \varepsilon)t$$
$$= A_S(0, Mx) - (c - \varepsilon)Mx + A_S(Mx, t) - (c - \varepsilon)(t - Mx),$$

and observe that $A_S(Mx, t) \stackrel{d}{=} A_S(0, t - Mx)$ since the process $A_S(0, t)$ is stationary. Thus, for $\delta > 0$ sufficiently small,

$$\mathbb{P}\left\{\sup_{t \geq Mx}\{A_S(0, t) - (c - \varepsilon)t\} > x\right\}$$

$$= \mathbb{P}\left\{\sup_{t \geq Mx}\{A_S(0, Mx) - (c - \varepsilon)Mx\right.$$

$$\left. + A_S(Mx, t) - (c - \varepsilon)(t - Mx)\} > x\right\}$$

$$= \mathbb{P}\left\{A_S(0, Mx) - (c - \varepsilon)Mx\right.$$

$$+ \sup_{t \geq Mx} \{A_S(Mx, t) - (c - \varepsilon)(t - Mx)\} > x\Big\}$$

$$\leq \mathbb{P}\{A_S(0, Mx) - (c - \varepsilon)Mx > -\delta(c - \varepsilon)Mx\}$$

$$+ \mathbb{P}\Big\{\sup_{t \geq Mx} \{A_S(0, t - Mx) - (c - \varepsilon)(t - Mx)\} > (1 + \delta(c - \varepsilon)M)x\Big\}$$

$$= \mathbb{P}\{A_S(0, Mx) > (1 - \delta)(c - \varepsilon)Mx\}$$

$$+ \mathbb{P}\Big\{\sup_{t \geq Mx} \{A_S(0, t - Mx) - (c - \varepsilon)(t - Mx)\} > (1 + \delta(c - \varepsilon)M)x\Big\}$$

$$\leq \mathbb{P}\Big\{\sup_{t \geq 0} \{A_S(0, t) - (1 - 2\delta)(c - \varepsilon)t\} > \delta(c - \varepsilon)Mx\Big\}$$

$$+ \mathbb{P}\Big\{\sup_{t \geq 0} \{A_S(0, t) - (c - \varepsilon)t\} > (1 + \delta(c - \varepsilon)M)x\Big\}$$

$$= \mathbb{P}\{\mathbf{V}_S^{(1 - 2\delta)(c - \varepsilon)} > \delta(c - \varepsilon)Mx\} + \mathbb{P}\{\mathbf{V}_S^{c - \varepsilon} > (1 + \delta(c - \varepsilon)M)x\}.$$

Using Corollary 4.2, for $\delta > 0$ sufficiently small,

$$\frac{\mathbb{P}\{\sup_{t \geq Mx}\{A_S(0, t) - (c - \varepsilon)t\} > x\}}{\mathbb{P}\{\mathbf{V}_S^c > x\}}$$

$$\leq \frac{K^{(2)}P_S(\delta(c - \varepsilon)Mx)}{K^{(1)}P_S(x)} + \frac{K^{(2)}P_S((1 + \delta(c - \varepsilon)M)x)}{K^{(1)}P_S(x)}.$$

Now let $x \to \infty$ and then $M \to \infty$ [use the fact that $P_S(\cdot)$ is of regular variation].                                                                    □

**5. Reduced-load equivalence.** In this section we provide the proofs of the various reduced-load equivalence results stated in Section 3. The proofs of the complementing results for the reduced system are presented in Section 6. Section 5.1 treats the case of several weakly dominant sets, culminating in a proof of (3.4). See Theorem 5.1; this also gives the proof of the special case in Theorem 3.1. In Section 5.2 we extend the results to the case of additional instantaneous, heavy-tailed input.

5.1. *Proof of* (3.4). Recall that $\Upsilon$ denotes the collection of all weakly dominant sets, and that $\Omega$ represents the collection of all minimally critical sets.

ASSUMPTION 5.1. For any $y$ and $\delta > 0$,

$$F_S^c(\delta) := \liminf_{x \to \infty} \frac{\mathbb{P}\{\mathbf{V}_S^{c + \delta} > x + y\}}{\mathbb{P}\{\mathbf{V}_S^c > x\}}$$

is independent of $y$. In addition, $\lim_{\delta \downarrow 0} F_S^c(\delta) = 1$.

ASSUMPTION 5.2. *For any $y$ and $\delta > 0$,*

$$G_S^c(\delta) := \limsup_{x \to \infty} \frac{\mathbb{P}\{\mathbf{V}_S^{c-\delta} > x - y\}}{\mathbb{P}\{\mathbf{V}_S^c > x\}}$$

*is independent of $y$. In addition, $\lim_{\delta \downarrow 0} G_S^c(\delta) = 1$.*

ASSUMPTION 5.3. *For any $\varepsilon > 0$,*

$$\lim_{x \to \infty} \frac{\mathbb{P}\{\mathbf{V}_{\mathcal{I}_1}^{\rho_{\mathcal{I}_1} + \varepsilon} > x/\mathcal{N}\}}{\mathbb{P}\{\mathbf{V}_S^c > x\}} = 0.$$

ASSUMPTION 5.4. *For any $\varepsilon > 0$,*

$$H_S^c(\varepsilon) := \limsup_{x \to \infty} \frac{\prod_{j \in S} \mathbb{P}\{\mathbf{V}_j^{\rho_j + \varepsilon} > x/\mathcal{N}\}}{\mathbb{P}\{\mathbf{V}_S^c > x\}} < \infty.$$

ASSUMPTION 5.5. *For any pair of sets $S \in \Upsilon$, $E \in \Omega \setminus \Upsilon$, for any $\varepsilon > 0$,*

$$\lim_{x \to \infty} \frac{\prod_{j \in E} \mathbb{P}\{\mathbf{V}_j^{\rho_j + \varepsilon} > x/\mathcal{N}\}}{\mathbb{P}\{\mathbf{V}_S^c > x\}} = 0.$$

THEOREM 5.1 (Generalized reduced-load equivalence; weakly dominant sets). *Suppose the sets $S \in \Lambda$ satisfy Assumptions 5.1–5.5. Then*

$$\mathbb{P}\{\mathbf{V} > x\} \sim \sum_{S \in \Lambda} P\{\mathbf{V}_S^{c_S} > x\}.$$

PROOF. As before, the proof consists of a lower bound and an upper bound which asymptotically coincide. For compactness, denote $Q(x) := \sum_{S \in \Lambda} \mathbb{P}\{\mathbf{V}_S^{c_S} > x\}$.

(Lower bound.) From Lemma 4.3, for any $\delta > 0$ and $y \geq 0$,

$$\mathbb{P}\{\mathbf{V} > x\} \geq \sum_{S \in \Lambda} \mathbb{P}\{\mathbf{V}_S^{c_S + \delta} > x + y\} \mathbb{P}\{\mathbf{Z}_{\mathcal{I} \setminus S}^{\rho_{\mathcal{I}} - \delta} \leq y\}$$

$$- \sum_{S_1, S_2 \in \Lambda, S_1 \neq S_2} \prod_{j \in S_1 \cup S_2} \mathbb{P}\{\mathbf{V}_j^{\rho_j + \varepsilon} > x/\mathcal{N}\}.$$

Note that if $S_1, S_2 \in \Lambda$, $S_1 \neq S_2$, then $S_1 \cup S_2$ cannot be a minimally critical set, so that $S_1 \cup S_2 \notin \Lambda$.

Thus, using Assumptions 5.1 and 5.4, and the inequality

$$\frac{\sum_i a_i}{\sum_i b_i} \geq \min_i \frac{a_i}{b_i}$$

for $a_i, b_i > 0$, we obtain

$$\liminf_{x \to \infty} \frac{\mathbb{P}\{\mathbf{V} > x\}}{Q(x)}$$

$$\geq \liminf_{x \to \infty} \sum_{S \in \Lambda} \mathbb{P}\{\mathbf{Z}_{\mathcal{I} \setminus S}^{\rho_\mathcal{I} \setminus S - \delta} \leq y\} \frac{\mathbb{P}\{\mathbf{V}_S^{c_S + \delta} > x + y\}}{Q(x)}$$

$$- \sum_{S_1, S_2 \in \Lambda, S_1 \neq S_2} \limsup_{x \to \infty} \frac{\prod_{j \in S_1 \cup S_2} \mathbb{P}\{\mathbf{V}_j^{\rho_j + \varepsilon} > x / \mathcal{N}\}}{Q(x)}$$

$$\geq \liminf_{x \to \infty} \min_{S \in \Lambda} \mathbb{P}\{\mathbf{Z}_{\mathcal{I} \setminus S}^{\rho_\mathcal{I} \setminus S - \delta} \leq y\} \frac{\mathbb{P}\{\mathbf{V}_S^{c_S + \delta} > x + y\}}{\mathbb{P}\{\mathbf{V}_S^{c_S} > x\}}$$

$$\geq \min_{S \in \Lambda} \mathbb{P}\{\mathbf{Z}_{\mathcal{I} \setminus S}^{\rho_\mathcal{I} \setminus S - \delta} \leq y\} \liminf_{x \to \infty} \frac{\mathbb{P}\{\mathbf{V}_S^{c_S + \delta} > x + y\}}{\mathbb{P}\{\mathbf{V}_S^{c_S} > x\}}$$

$$= \min_{S \in \Lambda} F_S^{c_S}(\delta) \mathbb{P}\{\mathbf{Z}_{\mathcal{I} \setminus S}^{\rho_\mathcal{I} \setminus S - \delta} \leq y\}.$$

Letting $y \to \infty$, then $\delta \downarrow 0$, we obtain

$$\liminf_{x \to \infty} \frac{\mathbb{P}\{\mathbf{V} > x\}}{Q(x)} \geq 1,$$

which completes the proof of the lower bound.

(Upper bound.)  From Lemma 4.4, for any $\delta > 0$ and $y$,

$$\mathbb{P}\{\mathbf{V} > x\} \leq \sum_{S \in \Lambda} \mathbb{P}\{\mathbf{V}_S^{c_S - \delta} > x - y\} + \mathbb{P}\{\mathbf{V}_{\mathcal{I}_1}^{\rho_{\mathcal{I}_1} + \varepsilon} > x / \mathcal{N}\}$$

$$+ \sum_{S \in \Lambda} \mathbb{P}\{\mathbf{V}_{\mathcal{I} \setminus S}^{\rho_\mathcal{I} \setminus S + \delta} > y\} \prod_{j \in S} \mathbb{P}\{\mathbf{V}_j^{\rho_j + \varepsilon} > x / \mathcal{N}\}$$

$$+ \sum_{E \in \Omega \setminus \Lambda} \prod_{j \in E} \mathbb{P}\{\mathbf{V}_j^{\rho_j + \varepsilon} > x / \mathcal{N}\}.$$

Thus, using Assumptions 5.2–5.5, and the inequality

$$\frac{\sum_i a_i}{\sum_i b_i} \leq \max_i \frac{a_i}{b_i}$$

for $a_i, b_i > 0$,

$$\mathbb{P}\{\mathbf{V} > x\} \leq \limsup_{x \to \infty} \sum_{S \in \Lambda} \frac{\mathbb{P}\{\mathbf{V}_S^{c_S - \delta} > x - y\}}{Q(x)} + \limsup_{x \to \infty} \frac{\mathbb{P}\{\mathbf{V}_{\mathcal{I}_1}^{\rho_{\mathcal{I}_1} + \varepsilon} > x / \mathcal{N}\}}{Q(x)}$$

$$+ \sum_{S \in \Lambda} \mathbb{P}\{\mathbf{V}_{\mathcal{I} \setminus S}^{\rho_\mathcal{I} \setminus S + \delta} > y\} \limsup_{x \to \infty} \frac{\prod_{j \in S} \mathbb{P}\{\mathbf{V}_j^{\rho_j + \varepsilon} > x / \mathcal{N}\}}{Q(x)}$$

$$+ \sum_{E \in \Omega \setminus \Lambda} \limsup_{x \to \infty} \frac{\prod_{j \in E} \mathbb{P}\{\mathbf{V}_j^{\rho_j + \varepsilon} > x/\mathcal{N}\}}{Q(x)}$$

$$\leq \limsup_{x \to \infty} \max_{S \in \Lambda} \frac{\mathbb{P}\{\mathbf{V}_S^{c_S - \delta} > x - y\}}{\mathbb{P}\{\mathbf{V}_S^{c_S} > x\}}$$

$$+ \sum_{S \in \Lambda} \mathbb{P}\{\mathbf{V}_{\mathcal{I} \setminus S}^{\rho_{\mathcal{I} \setminus S} + \delta} > y\} \limsup_{x \to \infty} \frac{\prod_{j \in S} \mathbb{P}\{\mathbf{V}_j^{\rho_j + \varepsilon} > x/\mathcal{N}\}}{\mathbb{P}\{\mathbf{V}_S^{c_S} > x\}}$$

$$\leq \max_{S \in \Lambda} \limsup_{x \to \infty} \frac{\mathbb{P}\{\mathbf{V}_S^{c_S - \delta} > x - y\}}{\mathbb{P}\{\mathbf{V}_S^{c_S} > x\}} + \sum_{S \in \Lambda} H_S(\varepsilon) \mathbb{P}\{\mathbf{V}_{\mathcal{I} \setminus S}^{\rho_{\mathcal{I} \setminus S} + \delta} > y\}$$

$$= \max_{S \in \Lambda} G_S^{c_S}(\delta) + \sum_{S \in \Lambda} H_S(\varepsilon) \mathbb{P}\{\mathbf{V}_{\mathcal{I} \setminus S}^{\rho_{\mathcal{I} \setminus S} + \delta} > y\}.$$

Letting $y \to \infty$, then $\delta \downarrow 0$, we obtain

$$\limsup_{x \to \infty} \frac{\mathbb{P}\{\mathbf{V} > x\}}{Q(x)} \leq 1,$$

which completes the proof. $\quad \square$

In order to complete the proof of the reduced-load equivalence result (3.1), it remains to be shown that a dominant set $S^* \subseteq \mathcal{I}_2$ with $A_j(\cdot) \in \mathcal{R}$ for all $j \in S^*$ satisfies Assumptions 5.1–5.5. That is done in the following two propositions for $S = S^*$.

PROPOSITION 5.1. *Let $S \subseteq \mathcal{I}_2$. If $A_j(\cdot) \in \mathcal{R}$ for all $j \in S$, then Assumptions 5.1 and 5.2 are satisfied for any $c \in (r_S - \min_{j \in S}\{r_j - \rho_j\}, r_S)$.*

PROOF. We first prove that Assumption 5.2 is satisfied. It follows from Theorem 6.4 (see also Corollary 6.1; it is important to note here that results from Section 6 do not rely on the results of this section) that if $A_j(\cdot) \in \mathcal{R}$ for all $j \in S$, then $\mathbb{P}\{\mathbf{V}_S^c > x\} \in \mathcal{IR}$. Since $\mathcal{IR} \subset \mathcal{L}$, it suffices to prove that the assumption is satisfied for $y = 0$.

Let $\varepsilon \in [0, r_S - c)$, and let $\delta \in (0, \varepsilon]$. Then

$$\mathbb{P}\{\mathbf{V}_S^{c-\delta} > x\} = \mathbb{P}\Big\{\sup_{t \geq 0}\{A_S(0,t) - (c - \delta)t\} > x\Big\}$$

$$\leq \mathbb{P}\Big\{\sup_{t \leq x\delta^{-1/2}} \{A_S(0,t) - (c - \delta)t\} > x\Big\}$$

$$+ \mathbb{P}\Big\{\sup_{t \geq x\delta^{-1/2}} \{A_S(0,t) - (c - \delta)t\} > x\Big\}$$

$$\leq \mathbb{P}\left\{\sup_{t \leq x\delta^{-1/2}} \{A_S(0,t) - ct\} > (1 - \delta^{1/2})x\right\}$$

$$+ \mathbb{P}\left\{\sup_{t \geq x\delta^{-1/2}} \{A_S(0,t) - (c - \varepsilon)t\} > x\right\}.$$

Thus,

$$\limsup_{x \to \infty} \frac{\mathbb{P}\{\mathbf{V}_S^{c-\delta} > x\}}{\mathbb{P}\{\mathbf{V}_S^c > x\}} \leq \limsup_{x \to \infty} \frac{\mathbb{P}\{\mathbf{V}_S^c > (1 - \delta^{1/2})x\}}{\mathbb{P}\{\mathbf{V}_S^c > x\}}$$

$$+ \limsup_{x \to \infty} \frac{\mathbb{P}\{\sup_{t \geq x\delta^{-1/2}} \{A_S(0,t) - (c - \varepsilon)t\} > x\}}{\mathbb{P}\{\mathbf{V}_S^c > x\}}.$$

The fact that $\mathbb{P}\{\mathbf{V}_S^c > x\} \in \mathcal{IR}$ implies that the first term tends to 1 as $\delta \downarrow 0$, while Lemma 4.5 (with $M = \delta^{-1/2}$) shows that the second term then goes to 0.

The proof that Assumption 5.1 holds is similar, and therefore omitted. $\square$

PROPOSITION 5.2. *Let $S \subseteq \mathcal{I}_2$. If $A_j(\cdot) \in \mathcal{R}$ for all $j \in S$, then Assumptions 5.3 and 5.4 are satisfied for any $c > \rho_S$. If, in addition, $S$ is a weakly dominant set, then Assumption 5.5 is satisfied as well.*

PROOF. Using Lemma 4.1,

$$\mathbb{P}\{\mathbf{V}_S^c > x\} \geq \prod_{j \in S} p_j \mathbb{P}\left\{\mathbf{A}_j^r > \frac{x}{r_S - c}\right\}.$$

Assumption 5.3 then follows from combining Assumption 2.1 and the assumption that $A_j(\cdot) \in \mathcal{R}$ for all $j \in S$.

Theorem 2.1 gives

$$\mathbb{P}\{\mathbf{V}_j^{\rho_j + \varepsilon} > x/\mathcal{N}\} \sim (1 - p_j)\frac{\rho_j}{\varepsilon}\mathbb{P}\left\{\mathbf{A}_j^r > \frac{x/\mathcal{N}}{r_j - \rho_j - \varepsilon}\right\}$$

for all $j \in \mathcal{I}_2$.

Assumption 5.4 then follows from the assumption that $A_j(\cdot) \in \mathcal{R}$ for all $j \in S$, and so does Assumption 5.5 in case $S$ is a weakly dominant set. $\square$

5.2. *Additional instantaneous input.* So far we have considered a scenario with only *fluid* heavy-tailed input. We now extend the reduced-load equivalence to the case with additional *instantaneous*, heavy-tailed input. We thus allow for an additional subset of flows $\mathcal{I}_3 \subseteq \mathcal{I}$ which generate instantaneous traffic bursts according to independent renewal processes. The interarrival times between bursts of flow $i$ are generally distributed with mean $1/\lambda_i$. The burst sizes $\mathbf{B}_i$ have a heavy-tailed distribution $B_i(\cdot)$ with mean $\beta_i < \infty$. Thus the traffic intensity of flow $i$ is $\rho_i := \lambda_i \beta_i$.

For each flow $i \in \mathcal{I}_3$, we assume that the burst size distribution is regularly varying of index $-\nu_i$, that is, $B_i(\cdot) \in \mathcal{R}_{-\nu_i}$ for some $\nu_i > 1$. The next result which is due to [32] then gives the tail behavior of the workload distribution for a single flow $i \in \mathcal{I}_3$ served in isolation.

THEOREM 5.2. *If $B_i^r(\cdot) \in \mathcal{S}$, $\rho_i < c$, then*

$$\tag{5.1} \mathbb{P}\{\mathbf{V}_i^c > x\} \sim \frac{\rho_i}{c - \rho_i} \mathbb{P}\{\mathbf{B}_i^r > x\}.$$

In order to formulate the results, we need to extend the concept of dominance introduced in Section 3. A flow $i \in \mathcal{I}_3$ is said to (weakly) dominate a flow $j \in \mathcal{I}_3$ if $\nu_i < (\leq) \nu_j$. A flow $i \in \mathcal{I}_3$ is said to (weakly) dominate a critical set $S \subseteq \mathcal{I}_2$ if $\nu_i - 1 < (\leq) \sum_{j \in S}(\nu_j - 1)$. A critical set $S \subseteq \mathcal{I}_2$ is said to (weakly) dominate a flow $i \in \mathcal{I}_3$ if $\nu_i - 1 > (\geq) \sum_{j \in S}(\nu_j - 1)$.

A flow $i \in \mathcal{I}_3$ is called (weakly) dominant if it (weakly) dominates all other flows $j \in \mathcal{I}_3$ as well as all critical sets $S \subseteq \mathcal{I}_2$. A critical set $S \subseteq \mathcal{I}_2$ is called (weakly) dominant if it (weakly) dominates any other critical set $U \subseteq \mathcal{I}_2$ as well as all flows $j \in \mathcal{I}_3$.

THEOREM 5.3. *Let $\mathcal{K} \subseteq \mathcal{I}_3$ and $\Upsilon \subseteq 2^{\mathcal{I}_2}$ be the collection of all weakly dominant flows and all weakly dominant sets, respectively. If $B_i(\cdot) \in \mathcal{R}$ for all $i \in \mathcal{K}$, and $A_j(\cdot) \in \mathcal{R}$ for all $j \in S$, $S \in \Upsilon$, then*

$$\tag{5.2} \mathbb{P}\{\mathbf{V} > x\} \sim \sum_{i \in \mathcal{K}} \mathbb{P}\{\mathbf{V}_i^{c_i} > x\} + \sum_{S \in \Upsilon} \mathbb{P}\{\mathbf{V}_S^{c_S} > x\},$$

*with $\mathbb{P}\{\mathbf{V}_i^{c_i} > x\}$ and $\mathbb{P}\{\mathbf{V}_S^{c_S} > x\}$ as in (5.1) and (3.2), (3.3), respectively.*

The proof of the above theorem is similar to that of Theorem 5.1 after a few modifications to Lemmas 4.3 and 4.4.

It may be worth mentioning that Theorem 5.3 continues to hold under the condition $B_i^r(\cdot) \in \mathcal{S}$ for all $i \in \mathcal{K}$, provided there are no weakly dominant sets of On–Off flows (the concept of dominance may be extended to subexponential distributions in a straightforward way). In particular, when there are simply no On–Off flows at all, one obtains the extension of Theorem 5.2 to the single-server queue fed by a superposition of renewal processes (which is not a renewal process). This result was obtained as Theorem 4.1 in [4], using a different approach.

Theorem 5.3 also provides an extension of a recent result in [11], who study an $M/G/1$ queue with two different speeds of service using an analytic approach. A queue with two service speeds fits into our framework by the observation that the varying service capacity can be regulated by an On–Off source.

**6. Tail asymptotics for the reduced system.** In this section we derive the tail asymptotics for the reduced system. In particular, we give a proof of (3.2) and (3.3).

For notational convenience, let $c$ be the capacity of the reduced system, let the set of flows be indexed as $\mathcal{J} = \{1, \ldots, N\}$, and denote $r := r_{\mathcal{J}}$ and $A(0,t) := A_{\mathcal{J}}(0,t)$. By definition, the reduced system satisfies the following two properties:

(i) The On period distribution of flow $i$ is regularly varying of index $-\nu_i < -1$, that is, $A_i(\cdot) \in \mathcal{R}_{-\nu_i}$.

(ii) All flows must be On for the drift of the workload process to be positive, that is, $c \in (r - \min_{i=1,\ldots,N} \{r_i - \rho_i\}, r)$.

We now state our main theorem.

THEOREM 6.1. *Consider a queue of capacity $c$ fed by $N$ On–Off flows. If $c \in (r - \min_{i=1,\ldots,N}\{r_i - \rho_i\}, r)$ with $r = \sum_{i=1}^{N} r_i$, and $A_j(\cdot) \in \mathcal{R}$ for all $j = 1, \ldots, N$, then*

$$\mathbb{P}\{\mathbf{V}^c > x\} \sim \left(\prod_{j=1}^{N} p_j\right) \sum_{\mathcal{J}_0 \subseteq \{1, \ldots, N\}} P_{\mathcal{J}_0}(x),$$

*where $P_{\mathcal{J}_0}(x)$ is given by (with $\mathcal{J}_1 = \{1, \ldots, N\} \setminus \mathcal{J}_0$)*

$$
\begin{aligned}
P_{\mathcal{J}_0}(x) = {} & \frac{1}{\prod_{i \in \mathcal{J}_1} \mathbb{E}\{\mathbf{A}_i\}} \\
& \times \int_{y_i \in (0,\infty), i \in \mathcal{J}_1} \prod_{i \in \mathcal{J}_1} \mathbb{P}\left\{(r-c)\mathbf{A}_i > \sum_{j \in \mathcal{J}_1} y_j(r_j - \rho_j) - (r-c)y_i + x\right\} \\
& \times \prod_{i \in \mathcal{J}_0} \mathbb{P}\left\{(r-c)\mathbf{A}_i^r > \sum_{j \in \mathcal{J}_1} y_j(r_j - \rho_j) + x\right\} \prod_{i \in \mathcal{J}_1} dy_i.
\end{aligned}
$$

(6.1)

An asymptotic characterization of $P_{\mathcal{J}_0}(x)$ which may be useful for further analysis is provided in Section 6.4. This characterization also shows that $\mathbb{P}\{\mathbf{V}^c > x\}$ and $P_{\mathcal{J}_0}(x)$ are regularly varying, and gives an expression for the pre-factor in the asymptotic expansion of $\mathbb{P}\{\mathbf{V}^c > x\}$.

The remainder of this section is organized as follows. Detailed heuristic arguments are given in Section 6.1. In Section 6.2, we prove some preliminary results on the most probable behavior of the process $\{A(0,t) - ct\}$. The proof of Theorem 6.1 is then completed in Section 6.3. Section 6.4 deals with the asymptotic behavior of $P_{\mathcal{J}_0}(x)$.

6.1. *Heuristic arguments.* The proof of Theorem 6.1 is quite lengthy. Nevertheless, it is based on a simple intuitive argument: the most likely way for $\mathbf{V}^c \equiv \sup_{t \geq 0}\{A(0, t) - ct\}$ to reach a large value is that all flows have been simultaneously On for a long time. Specifically, each flow is likely to contribute through *exactly one* "long" On period; apart from these long On periods, the flows show typical behavior.

The above heuristic argument may be used for computing $\sup_{t \geq 0}\{A(0, t) - ct\}$. Let us say that the long On period of flow $i$ begins at time $s_i$ and ends at time $s_i + t_i$. Define

$$t^* := \min_{i=1,\ldots,N}\{s_i + t_i\}$$

as the time epoch at which the first of the long On periods finishes. One may also interpret $t^*$ as the time epoch at which the process $\{A(0, t) - ct\}$ reaches its largest value. Note that $A_i(0, s_i) \approx \rho_i s_i$, $A_i(s_i, s_i + t_i) = r_i t_i$, and $A_i(s_i + t_i, s_i + t_i + t) \approx \rho_i t$, $t \geq 0$. One thus obtains, using the fact that $c \in (r - \min_{i=1,\ldots,N}\{r_i - \rho_i\}, r)$,

$$
\begin{aligned}
\sup_{t \geq 0}\{A(0, t) - ct\} &\approx A(0, t^*) - ct^* \\
&\approx \sum_{i=1}^{N}[\rho_i s_i + r_i(t^* - s_i)] - ct^* \\
&= \sum_{i=1}^{N}(\rho_i - r_i)s_i + (r - c)t^*.
\end{aligned}
$$

(6.2)

The problem is thus reduced to calculating

$$(6.3) \qquad \mathbb{P}\left\{\sum_{i=1}^{N}(\rho_i - r_i)s_i + (r - c)\min_{i=1,\ldots,N}\{s_i + t_i\} > x\right\}.$$

Although the proof is based on the representation $\mathbf{V}^c \equiv \sup_{t \geq 0}\{A(0, t) - ct\}$, it is useful to keep the original workload process $\sup_{0 \leq s \leq t}\{A(s, t) - c(t - s)\}$ in mind as well. Figure 1 shows a typical scenario leading to a large workload level (so small fluctuations are ignored) in the case of two On–Off flows.

At a certain time $\omega_0$, the first long On period begins. Before that time, both flows show average behavior. The queue starts to build (at rate $r_1 + r_2 - c$) at time $\omega_1$ when the second long On period begins, and reaches its largest level at time $\omega_3$. Level $x$ is crossed at time $\omega_2$.

Between times $\omega_3$ and $\omega_4$, the queue drains at rate $c - r_1 - \rho_2$: flow 1 is still in the middle of its long On period, and flow 2 shows average behavior (remember small fluctuations are neglected). The process is still above level $x$ between times $\omega_4$ and $\omega_5$. However, here both flows show average behavior again, causing a negative drift $c - \rho_1 - \rho_2$.

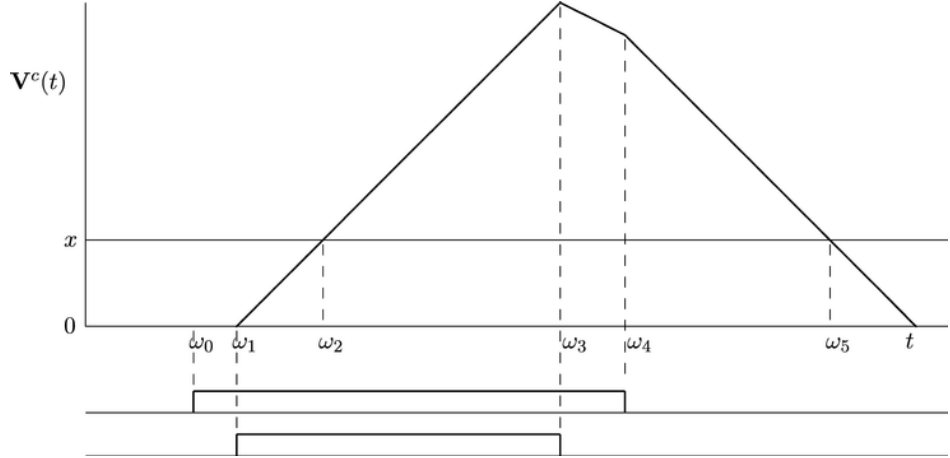

Fig. 1.  *Typical overflow scenario for two On–Off flows.*

The figure illustrates why the analysis of the reduced system is still quite complicated:

   (i)  Although the long On periods must significantly overlap, the difference between the finishing times of these On periods can be quite large (of order $x$, hence not negligible).

   (ii)  Given that the observed workload is larger than $x$, it is not necessarily the case that all flows are in the middle of their long On periods. In Figure 1, this is only the case in the time interval $(\omega_2, \omega_3)$. In fact, for any given flow, its long On period may have finished a long time ago. Consequently, there are $2^N$ different possibilities (corresponding to which subsets of the flows are still in the middle of their long On periods). Sample-path wise, there are $N + 1$ different time intervals in which the workload may be larger than $x$ (depending on how many of the flows are still in the middle of their long On periods).

   (iii)  Specifically, given that the observed workload is larger than $x$, it may still have been even larger at an earlier time epoch. In Figure 1, this is the case in the time intervals $(\omega_3, \omega_4)$ and $(\omega_4, \omega_5)$.

These complications do not arise if one considers a related problem, which concerns the overflow probability in a fluid queue with a *finite buffer* of size $x$. As is shown in [23], the analysis of the reduced system is then considerably simpler. It suffices to use bounds which are similar to Lemmas 4.1 and 4.2, and to combine these with the asymptotic results for a single On–Off flow in [21, 40].

6.2. *Characterization of most probable behavior.* In this section we prove some preliminary results characterizing the most probable behavior of the process $\{A(0,t) - ct\}$ given that it reaches a large value. In particular, we formalize the following two heuristic statements, resulting in a formal version of (6.2).

(i) Each flow contributes to $\sup_{t\geq 0}\{A(0,t) - ct\}$ through exactly one "long" On period.

(ii) Apart from these long On periods, the flows show typical behavior.

An On period is referred to as "long" when larger than $\varepsilon x$, with $\varepsilon$ some small but positive constant. In order to formalize the above statements, we need to keep track of how many such long On periods occur.

With that in mind, we define $\mathcal{N}_i(A, B)$, for intervals $A, B \subseteq [0, \infty)$, as the number of On periods of flow $i$ of which the length is contained in $A$ and which overlap (in time) with $B$. For compactness, denote

$$\mathcal{N}_i(u,t) \equiv \mathcal{N}_i((u,\infty), [0,t]).$$

We now proceed with a few preparatory lemmas.

First we show how to obtain an upper bound for the workload process in terms of a simple random walk. As in (4.1), we have $V^c(t) \leq V_i^{d_i}(t)$ for all $i = 1, \ldots, N$, with $d_i := c - r_{\mathcal{I}\setminus\{i\}} = c - r + r_i$. Recall that $V_i^{d_i}(t) \overset{d}{=} \sup_{0\leq s\leq t}\{A_i(0,s) - d_i s\}$. Now let, for fixed $i$, $\mathbf{S}_{in} = \mathbf{X}_{i1} + \cdots + \mathbf{X}_{in}$ be a random walk with step sizes $\mathbf{X}_{im} = (r_i - d_i)\mathbf{A}_{im} - d_i\mathbf{U}_{im}$, with $\mathbf{A}_{im}$ and $\mathbf{U}_{im}$ i.i.d. random variables distributed as the On and Off periods of flow $i$, respectively.

Since $c \in (r - \min_{i=1,\ldots,N}\{r_i - \rho_i\}, r)$, we have $\rho_i < d_i$ for all $i = 1, \ldots, N$, so that $\mathbb{E}\{\mathbf{X}_{i1}\} < 0$, that is, the random walk has negative drift. Because of the sawtooth nature of the process $A_i(0, s) - d_i s$, we have

$$\sup_{0\leq s\leq t}\{A_i(0,s) - d_i s\} \leq (r_i - d_i)(\mathbf{B}_i\mathbf{A}_{i0}^r + (1 - \mathbf{B}_i)\mathbf{A}_{i0}) + \sup_{n\leq N_i^A(t)}\mathbf{S}_{in},$$

with $N_i^A(t)$ denoting the number of Off periods of flow $i$ elapsed during $[0, t]$ which started after time 0 [for a formal definition see (4.5)].

The above observations are summarized in the following auxiliary lemma.

LEMMA 6.1. *For all $\varepsilon > 0$, $t$ and $x$,*

$$\mathbb{P}\{V^c(t) > x, \mathcal{N}_i(\varepsilon x, t) = 0\}$$

$$\leq \mathbb{P}\left\{\sup_{n\leq N_i^A(t)}\mathbf{S}_{in} > x(1 - \varepsilon(r_i - d_i)), \mathcal{N}_i(\varepsilon x, t) = 0\right\}.$$

Proof.    We have

$$\mathbb{P}\{V^c(t) > x, \mathcal{N}_i(\varepsilon x, t) = 0\}$$

$$\leq \mathbb{P}\{V_i^{d_i}(t) > x, \mathcal{N}_i(\varepsilon x, t) = 0\}$$

$$\leq \mathbb{P}\bigg\{(r_i - d_i)(\mathbf{B}_i \mathbf{A}_{i0}^r + (1 - \mathbf{B}_i)\mathbf{A}_{i0})$$

$$+ \sup_{n \leq N_i^A(t)} \mathbf{S}_{in} > x, \mathcal{N}_i(\varepsilon x, t) = 0\bigg\}$$

$$\leq \mathbb{P}\bigg\{\sup_{n \leq N_i^A(t)} \mathbf{S}_{in} > x(1 - \varepsilon(r_i - d_i)), \mathcal{N}_i(\varepsilon x, t) = 0\bigg\}.$$

The last inequality follows from the fact that $\mathbf{A}_{i0}^r$ and $\mathbf{A}_{i0}$ must be smaller than $\varepsilon x$ if $\mathcal{N}_i(\varepsilon x, t) = 0$. $\quad\square$

To obtain upper bounds for probabilities as in Lemma 6.1, we will frequently apply the following key lemma, which is a trivial modification of Lemma 3 in [34].

Lemma 6.2.    Let $\mathbf{S}_n = \mathbf{X}_1 + \cdots + \mathbf{X}_n$ be a random walk with i.i.d. step sizes such that $\mathbb{E}\{\mathbf{X}_1\} < 0$ and $\mathbb{E}\{(\mathbf{X}_1^+)^p\} < \infty$ for some $p > 1$. Then, for any $\beta < \infty$, there exists an $\varepsilon^* > 0$ and a function $\phi(\cdot) \in \mathcal{R}_{-\beta}$ such that, for $\varepsilon \in (0, \varepsilon^*]$,

$$\mathbb{P}\{\mathbf{S}_n > x \mid \mathbf{X}_j \leq \varepsilon x, j = 1, \ldots, n\} \leq \phi(x),$$

for all $n$ and all $x$.

Note that if $\mathbf{X}_j$ can be represented as the difference of two nonnegative independent random variables $\mathbf{X}_j^1$ and $\mathbf{X}_j^2$, then the lemma remains valid if the $\mathbf{X}_j$'s are replaced by $\mathbf{X}_j^1$.

The final preparatory lemma is a simple consequence of Corollary 4.2, which will be used several times in combination with Lemma 6.2 that probabilities of certain events are of $o(\mathbb{P}\{\mathbf{V}^c > x\})$. Define

$$P(x) := \prod_{j=1}^N \mathbb{P}\{\mathbf{A}_j^r > x\} \in \mathcal{R}_{-\mu}, \qquad \mu := \sum_{j=1}^N (\nu_j - 1).$$

Lemma 6.3.

$$\limsup_{x \to \infty} \frac{P(x)}{\mathbb{P}\{\mathbf{V}^c > x\}} < \infty.$$

We now show that, with overwhelming probability (as $x \to \infty$), the rare event $\{\mathbf{V}^c > x\}$ occurs as follows:

(i) The process $\{A(0,t) - ct\}$ reaches level $x$ before time $Mx$ for some large $M$.

(ii) Up to time $Mx$, each flow generates *exactly one* long On period, that is, $\mathcal{N}_i(\varepsilon x, t) = 1$ for $i = 1, \ldots, N$.

PROPOSITION 6.1.

$$\lim_{M \to \infty} \liminf_{x \to \infty} \frac{\mathbb{P}\{V^c(Mx) > x\}}{\mathbb{P}\{\mathbf{V}^c > x\}} = 1.$$

PROOF. By definition,

$$\mathbb{P}\{\mathbf{V}^c > x\} = \mathbb{P}\Big\{\sup_{t \geq 0}\{A(0,t) - ct\} > x\Big\}$$

$$\leq \mathbb{P}\Big\{\sup_{0 \leq t \leq Mx}\{A(0,t) - ct\} > x\Big\} + \mathbb{P}\Big\{\sup_{t \geq Mx}\{A(0,t) - ct\} > x\Big\}$$

$$= \mathbb{P}\{V^c(Mx) > x\} + \mathbb{P}\Big\{\sup_{t \geq Mx}\{A(0,t) - ct\} > x\Big\}.$$

Thus, it suffices to show

$$\lim_{M \to \infty} \limsup_{x \to \infty} \frac{\mathbb{P}\{\sup_{t \geq Mx}\{A(0,t) - ct\} > x\}}{\mathbb{P}\{\mathbf{V}^c > x\}} = 0,$$

which, however, follows directly from Lemma 4.5.   □

Now suppose that the workload reaches level $x$. By the previous proposition, we may assume that this occurs before time $Mx$ (for $M$ sufficiently large). The next two propositions show that we may restrict the attention to a scenario where each flow initiates *exactly one* long On period before time $Mx$.

The first proposition indicates that each flow has *at least* one long On period.

PROPOSITION 6.2. *For all $i$ and all $M$, there exists an $\varepsilon^* > 0$ such that, for all $\varepsilon \in (0, \varepsilon^*]$,*

$$\mathbb{P}\{V^c(Mx) > x, \mathcal{N}_i(\varepsilon x, Mx) = 0\} = o(\mathbb{P}\{\mathbf{V}^c > x\}),$$

*as $x \to \infty$.*

PROOF.   Define $N_i^U(t) := \max\{n : \sum_{j=1}^n \mathbf{U}_{ij} \leq t\} + 1$. Note that $N_i^A(t) \leq N_i^U(t)$.

Using Lemma 6.1, taking $t = Mx$,

$$\mathbb{P}\{\mathbf{V}^c(Mx) > x, \mathcal{N}_i(\varepsilon x, Mx) = 0\}$$

$$\leq \mathbb{P}\left\{ \sup_{n \leq N_i^A(Mx)} \mathbf{S}_{in} > x(1 - \varepsilon(r_i - d_i)), \mathcal{N}_i(\varepsilon x, Mx) = 0 \right\}$$

$$\leq \mathbb{P}\left\{ \sup_{n \leq N_i^A(Mx)} \mathbf{S}_{in} > x(1 - \varepsilon(r_i - d_i)) \,|\, \mathcal{N}_i(\varepsilon x, Mx) = 0 \right\}$$

$$= \mathbb{P}\left\{ \sup_{n \leq N_i^A(Mx)} \mathbf{S}_{in} > x(1 - \varepsilon(r_i - d_i)) \,|\, \mathbf{A}_{ij} < \varepsilon x, \right.$$
$$\left. j = 1, \ldots, N_i^A(Mx) \right\}$$

$$= \mathbb{P}\left\{ \sup_{n \leq N_i^A(Mx)} \mathbf{S}_{in} > x(1 - \varepsilon(r_i - d_i)) \,|\, \mathbf{A}_{ij} < \varepsilon x, j \geq 1 \right\}$$

$$\leq \mathbb{P}\left\{ \sup_{n \leq N_i^U(Mx)} \mathbf{S}_{in} > x(1 - \varepsilon(r_i - d_i)) \,|\, \mathbf{A}_{ij} < \varepsilon x, j \geq 1 \right\}$$

$$= \mathbb{P}\left\{ \sup_{n \leq N_i^U(Mx)} \mathbf{S}_{in} > x(1 - \varepsilon(r_i - d_i)) \,|\, \mathbf{A}_{ij} < \varepsilon x, \right.$$
$$\left. j = 1, \ldots, N_i^U(Mx) \right\}$$

$$\leq \mathbb{P}\left\{ \sup_{n \leq M_2 x} \mathbf{S}_{in} > x(1 - \varepsilon(r_i - d_i)) \,|\, \mathbf{A}_{ij} < \varepsilon x, j \geq 1 \right\}$$
$$+ \mathbb{P}\{N_i^U(Mx) > M_2 x\}.$$

The second term decays exponentially fast in $x$ if $M_2 > \lambda_i M$. The first term can be bounded by

$$\sum_{m=1}^{M_2 x} \mathbb{P}\{\mathbf{S}_{im} > x(1 - \varepsilon(r_i - d_i)) \,|\, \mathbf{A}_{ij} \leq \varepsilon x, j = 1, \ldots, m\}.$$

According to Lemma 6.2, there exists an $\varepsilon^* > 0$ and a function $\phi(\cdot) \in \mathcal{R}_{-\beta}$ with $\beta > \mu + 1$, such that, for $\varepsilon \in (0, \varepsilon^*]$, the last quantity is upper bounded by $M_2 x \phi(x)$. The latter function is regularly varying of index $1 - \beta < -\mu$. Invoking Lemma 6.3 then completes the proof.   □

The next proposition shows that each flow has *at most* one long On period.

PROPOSITION 6.3. *For all $i$, all $M$ and all $\varepsilon > 0$,*

$$\mathbb{P}\{V^c(Mx) > x, \mathcal{N}_i(\varepsilon x, Mx) \geq 2\} = o(\mathbb{P}\{\mathbf{V}^c > x\}),$$

*as $x \to \infty$.*

PROOF. Without loss of generality we may take $i = 1$. By Proposition 6.2 it suffices to show that

$$\mathbb{P}\{V^c(Mx) > x, \mathcal{N}_1(\varepsilon x, Mx) \geq 2, \mathcal{N}_i(\varepsilon x, Mx) \geq 1, i \geq 2\} = o(\mathbb{P}\{\mathbf{V}^c > x\}).$$

Note that the left-hand side is bounded by

$$\mathbb{P}\{\mathcal{N}_1(\varepsilon x, Mx) \geq 2\} \prod_{i=2}^{N} \mathbb{P}\{\mathcal{N}_i(\varepsilon x, Mx) \geq 1\}.$$

Thus, invoking Lemma 6.3, it suffices to show that:

(i) $\mathbb{P}\{\mathcal{N}_i(\varepsilon x, Mx) \geq 1\}$ is bounded by a function which is regularly varying of index $1 - \nu_i$.

(ii) $\mathbb{P}\{\mathcal{N}_i(\varepsilon x, Mx) \geq 2\} = o(\mathbb{P}\{\mathcal{N}_i(\varepsilon x, Mx) \geq 1\})$.

We will prove both assertions for $i = 1$. For assertion (i), note that

$$\mathbb{P}\{\mathcal{N}_1(\varepsilon x, Mx) \geq 1\}$$
$$\leq p_1 \mathbb{P}\{\mathbf{A}_1^r \geq \varepsilon x\} + \mathbb{P}\{\#\{j \in \{1, \ldots, N_1^U(Mx)\} : \mathbf{A}_{1j} \geq \varepsilon x\} \geq 1\}.$$

The first term is in $\mathcal{R}_{1-\nu_1}$. By conditioning upon $N_1^U(Mx)$, the second term can be bounded by $\mathbb{E}\{N_1^U(Mx)\}\mathbb{P}\{\mathbf{A}_1 \geq \varepsilon x\}$, which is also regularly varying of index $1 - \nu_1$. To prove assertion (ii), note that

$$\mathbb{P}\{\mathcal{N}_1(\varepsilon x, Mx) \geq 2\}$$
$$\leq p_1 \mathbb{P}\{\mathbf{A}_1^r \geq \varepsilon x\} \mathbb{P}\{\mathcal{N}_1((\varepsilon x, \infty), (0, Mx]) \geq 1\}$$
$$+ \mathbb{P}\{\mathcal{N}_1((\varepsilon x, \infty), (0, Mx]) \geq 2\}.$$

Using $\mathbb{P}\{\mathcal{N}_1((\varepsilon x, \infty), (0, Mx]) \geq 1\} \leq \mathbb{P}\{\mathcal{N}_1(\varepsilon x, Mx) \geq 1\}$ and assertion (i), it follows that the first term is of $o(\mathbb{P}\{\mathcal{N}_1(\varepsilon x, Mx) \geq 1\})$. To bound the second term, condition (again) on $N_1^U(Mx)$. This yields

$$\mathbb{P}\{\mathcal{N}_1((\varepsilon x, \infty), (0, Mx)) \geq 2\} \leq \mathbb{E}\{N_1^U(Mx)^2\}\mathbb{P}\{\mathbf{A}_1 \geq \varepsilon x\}^2.$$

Finally, note that $\mathbb{E}\{N_1^U(Mx)^2\}$ is quadratic in $x$ for $x \to \infty$. $\square$

We have now shown that, with overwhelming probability, each flow contributes to a large value of $\sup_{t \geq 0}\{A(0, t) - ct\}$ through exactly one long On period. We thus proceed with the second statement (as indicated at the beginning of this section), implying that apart from these long On periods, the flows show typical behavior. In order to formalize that statement, we need to introduce some notation. Define

$$\tau(y) := \inf\{t \geq 0 : A(0,t) - ct = y\}$$

as the first time at which the process $\{A(0,t) - ct\}$ reaches level $y$.

For fixed $\varepsilon > 0$ and $x$, let $\tau_{s,i}(\varepsilon x)$ and $\tau_{f,i}(\varepsilon x)$ be the respective starting and finishing times of the first On period of flow $i$ exceeding length $\varepsilon x$. Denote

$$\tau_s(\varepsilon x) := \max_{i=1,\ldots,N} \tau_{s,i}(\varepsilon x)$$

and

$$\tau_f(\varepsilon x) := \min_{i=1,\ldots,N} \tau_{f,i}(\varepsilon x).$$

Note that all flows are in the middle of their long On periods between times $\tau_s(\varepsilon x)$ and $\tau_f(\varepsilon x)$. We will show that the fluctuations of the process $\{A(0,t) - ct\}$ away from the mean before time $\tau_s(\varepsilon x)$ and after time $\tau_f(\varepsilon x)$ can be neglected.

More formally, the next two propositions show that, given that the workload reaches level $x$ before time $Mx$, there exists for any small $\delta > 0$ an $\varepsilon_\delta$ such that, for all $\varepsilon \in (0, \varepsilon_\delta)$,

$$\tau_s(\varepsilon x) \leq \tau(\delta x) < \tau((1-\delta)x) \leq \tau_f(\varepsilon x).$$

Thus, the workload remains small up to time $\tau_s(\varepsilon x)$, and reaches a level close to $x$ before time $\tau_f(\varepsilon x)$, as depicted in Figure 2.

The first proposition indicates that it is most unlikely that the process $\{A(0,t) - ct\}$ reaches level $\delta x$ before time $\tau_s(\varepsilon x)$.

PROPOSITION 6.4.  *For any $\delta > 0$, there exists an $\varepsilon^* > 0$ such that, for all $\varepsilon \in (0, \varepsilon^*]$,*

$$\mathbb{P}\{\tau(\delta x) < \tau_s(\varepsilon x)\} = o(\mathbb{P}\{\mathbf{V}^c > x\}).$$

PROOF.   For compactness, denote $\tau_s \equiv \tau_s(\varepsilon x)$, $\tau_{s,i} \equiv \tau_{s,i}(\varepsilon x)$. Then

$$\mathbb{P}\{\tau(\delta x) < \tau_s\} = \mathbb{P}\{V^c(\tau_s) > \delta x\} \leq \sum_{i=1}^{N} \mathbb{P}\{V^c(\tau_{s,i}) > \delta x\}.$$

We bound each term in the last summation.

Define $N_i(\varepsilon x) := N_i^A(\tau_{s,i}^-)$ as the number of On periods initiated by flow $i$ before the first On period exceeding length $\varepsilon x$. Note that $N_i(\varepsilon x) + 1$ is geometrically distributed with parameter $\mathbb{P}\{\mathbf{A}_i > \varepsilon x\}$.

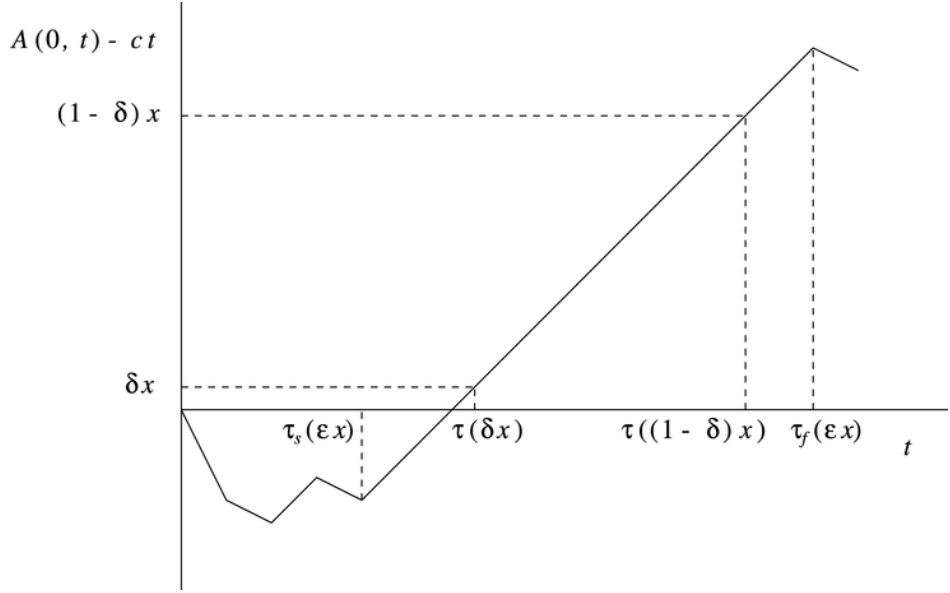

Fig. 2.   *Typical path to overflow.*

Using Lemma 6.1, taking $t = \tau_{s,i}$,

$$\mathbb{P}\{V^c(\tau_{s,i}) > \delta x\}$$

$$= \mathbb{P}\{V^c(\tau_{s,i}) > \delta x, \mathcal{N}_i((\varepsilon x, \infty), [0, \tau_{s,i})) = 0\}$$

$$\leq \mathbb{P}\left\{\sup_{n \leq N_i(\varepsilon x)} \mathbf{S}_n > x(\delta - \varepsilon(r_i - d_i)), \mathbf{A}_{ij} \leq \varepsilon x, j = 1, \ldots, N_i(\varepsilon x)\right\}$$

$$\leq \sum_{m=1}^{\infty} \mathbb{P}\{N_i(\varepsilon x) = m\}$$

$$\times \mathbb{P}\left\{\sup_{n \leq m} \mathbf{S}_n > x(\delta - \varepsilon(r_i - d_i)), \mathbf{A}_{ij} \leq \varepsilon x, j = 1, \ldots, m\right\}$$

$$\leq \sum_{m=1}^{\infty} \mathbb{P}\{N_i(\varepsilon x) = m\}$$

$$\times \mathbb{P}\left\{\sup_{n \leq m} \mathbf{S}_n > x(\delta - \varepsilon(r_i - d_i)) \mid \mathbf{A}_{ij} \leq \varepsilon x, j = 1, \ldots, m\right\}.$$

According to Lemma 6.2, there exists an $\varepsilon^* > 0$ and a function $\phi(\cdot) \in \mathcal{R}_{-\beta}$ with $\beta > 2\nu + 1$, such that, for $\varepsilon \in (0, \varepsilon^*)$, the last quantity is upper bounded by

$$\mathbb{E}\{N_i(\varepsilon x)\}\phi(x) = \frac{\phi(x)\mathbb{P}\{\mathbf{A}_i \leq \varepsilon x\}}{\mathbb{P}\{\mathbf{A}_i > \varepsilon x\}},$$

which is regularly varying of index $\nu_i - \beta < \mu + 1 - (2\mu + 1) = -\mu$.

Invoking Lemma 6.3 then completes the proof. □

The next proposition shows that, given that the process $\{A(0, t) - ct\}$ reaches level $x$ before time $Mx$, most probably level $(1 - \delta)x$ is crossed before time $\tau_f(\varepsilon x)$.

PROPOSITION 6.5. *For any $\delta > 0$ and $M < \infty$, there exists an $\varepsilon^* > 0$ such that, for all $\varepsilon \in (0, \varepsilon^*)$,*

$$\mathbb{P}\{\tau((1 - \delta)x) > \tau_f(\varepsilon x), V^c(Mx) > x\} = o(\mathbb{P}\{\mathbf{V}^c > x\}).$$

PROOF. For conciseness, denote $\tau_f \equiv \tau_f(\varepsilon x)$, $\tau_{f,i} \equiv \tau_{f,i}(\varepsilon x)$. By Propositions 6.2 and 6.3, it suffices to show that

$$\mathbb{P}\{\tau((1 - \delta)x) > \tau_f, V^c(Mx) > x, \mathcal{N}_i(\varepsilon x, Mx) = 1 \text{ for all } i = 1, \dots, N\}$$
$$= o(\mathbb{P}\{\mathbf{V}^c > x\}).$$

Note that

$$\mathbb{P}\{\tau((1 - \delta)x) > \tau_f, V^c(Mx) > x, \mathcal{N}_i(\varepsilon x, Mx) = 1 \text{ for all } i = 1, \dots, N\}$$
$$= \mathbb{P}\{V^c(\tau_f) > (1 - \delta)x, V^c(Mx) > x,$$
$$\mathcal{N}_i(\varepsilon x, Mx) = 1 \text{ for all } i = 1, \dots, N\}$$
$$\leq \sum_{i=1}^{N} \mathbb{P}\{V^c(\tau_{f,i}) > (1 - \delta)x, V^c(Mx) > x, \mathcal{N}_i(\varepsilon x, Mx) = 1\}.$$

As before, we bound each term in the last summation:

$$\mathbb{P}\{V^c(\tau_{f,i}) > (1 - \delta)x, V^c(Mx) > x, \mathcal{N}_i(\varepsilon x, Mx) = 1\}$$
$$\leq \mathbb{P}\Big\{\sup_{0 \leq t \leq \tau_{f,i}} \{A(0, t) - ct\} < (1 - \delta)x, \sup_{0 \leq t \leq Mx} \{A(0, t) - ct\} > x,$$
$$\mathcal{N}_i((\varepsilon x, \infty), (\tau_{f,i}, Mx]) = 0\Big\}$$
$$\leq \mathbb{P}\Big\{\sup_{\tau_{f,i} \leq t \leq Mx} \{A(\tau_{f,i}, t) - c(t - \tau_{f,i})\} > \delta x,$$
$$\mathcal{N}_i((\varepsilon x, \infty), (\tau_{f,i}, Mx]) = 0\Big\}$$
$$\leq \mathbb{P}\Big\{\sup_{\tau_{f,i} \leq t \leq Mx} \{A_i(\tau_{f,i}, t) - d_i(t - \tau_{f,i})\} > \delta x,$$
$$\mathcal{N}_i((\varepsilon x, \infty), (\tau_{f,i}, Mx]) = 0\Big\}.$$

The first inequality follows from the definitions. The second inequality follows from properties of the sup operator, and the last inequality is obtained by assuming that all flows but $i$ are On between times $\tau_{f,i}$ and $Mx$.

Note that the last probability is upper bounded by

$$\mathbb{P}\Bigg\{\sup_{N_i(\varepsilon x)+2 \leq n \leq N_i^A(Mx)} \mathbf{S}_n - \mathbf{S}_{N_i(\varepsilon x)+1} > \delta x,$$

$$\mathbf{A}_j \leq \varepsilon x, N_i(\varepsilon x) + 2 \leq j \leq N_i^A(Mx)\Bigg\}.$$

The latter probability can be upper bounded by a function which is regularly varying of index $-\beta < -\mu$ in a similar fashion as in the proof of Propositions 6.2 and 6.4.

The proof is complete by invoking Lemma 6.3. $\quad\square$

Propositions 6.4 and 6.5 may be used to obtain the following result.

COROLLARY 6.1. *If $A_j(\cdot) \in \mathcal{R}$ for all $j = 1, \dots, N$, then $\mathbb{P}\{\mathbf{V}^c > x\} \in \mathcal{IR}$.*

The above result suffices to prove the reduced-load equivalence. (See Section 5, in particular Proposition 5.1, for the details.) However, determining the exact asymptotic behavior of $\mathbb{P}\{\mathbf{V}^c > x\}$ requires further analysis, to be found in Sections 6.3 and 6.4. In particular, the analysis in Section 6.4 will lead to a sharper version of Corollary 6.1, showing that $\mathbb{P}\{\mathbf{V}^c > x\} \in \mathcal{R}$ (which is a strict subset of $\mathcal{IR}$).

Nevertheless, we sketch a direct proof of Corollary 6.1 which we believe is of independent interest. For the formal proof details we refer to [41].

SKETCH OF PROOF OF COROLLARY 6.1. The idea of the proof is as follows. If $\mathbf{V}^c > x$, then Propositions 6.4 and 6.5 show that the process $\{A(0,t) - ct\}$ reaches the level $(1-\delta)x$ after all flows have been On for at least $\frac{(1-2\delta)x}{r-c}$ time units. Since $A_j(\cdot) \in \mathcal{R} \subset \mathcal{IR}$ for all $j = 1, \dots, N$, with high probability, all flows remain On for at least $\frac{2\delta x}{r-c}$ more time units. This yields

$$\lim_{\delta\downarrow 0}\liminf_{x\to\infty}\mathbb{P}\{\mathbf{V}^c > (1+\delta)x \mid \mathbf{V}^c > x\} = 1,$$

implying the desired statement (by definition). $\quad\square$

6.3. *Proof of Theorem 6.1.* In this section we give a proof of Theorem 6.1. First we consolidate the key results from the previous section in the following theorem.

THEOREM 6.2. *For any $\delta > 0$, there exists an $\varepsilon^* > 0$ such that, for all $\varepsilon \in (0, \varepsilon^*)$,*

$$\mathbb{P}\{A(0, \tau_f(\varepsilon x)) - c\tau_f(\varepsilon x) > x\} \leq \mathbb{P}\{\mathbf{V}^c > x\}$$
$$\lesssim \mathbb{P}\{A(0, \tau_f(\varepsilon x)) - c\tau_f(\varepsilon x) > (1 - \delta)x\}.$$

PROOF. The lower bound is trivial. The upper bound follows from Propositions 6.1, 6.4 and 6.5.   □

In order to obtain tight bounds for the probabilities in Theorem 6.2, we condition upon $\tau_{s,i}$ for all $i$. Hence, for any $\mathcal{J}_0 \subseteq \mathcal{J}$, define the event $D_{\mathcal{J}_0}(\varepsilon x)$ by

$$D_{\mathcal{J}_0}(\varepsilon x) := \{\tau_{s,i}(\varepsilon x) = 0 \text{ for all } i \in \mathcal{J}_0; \tau_{s,i}(\varepsilon x) > 0 \text{ for all } i \notin \mathcal{J}_0\}.$$

The event $D_{\mathcal{J}_0}(\varepsilon x)$ implies that the flows $i \in \mathcal{J}_0$ started their long On period before time 0 (remember that we consider the system in stationarity). The flows $i \in \mathcal{J}_1$ start their long On period at a later time epoch.

Denote $\mathbb{P}_{\mathcal{J}_0}\{\cdot\} = \mathbb{P}\{\cdot \mid D_{\mathcal{J}_0}(\varepsilon x)\}$. The following two lemmas will be useful for providing tight upper and lower bounds for the probabilities in Theorem 6.2.

LEMMA 6.4 (Upper bound). *For any $\delta > 0$, there exists an $\varepsilon_\delta > 0$ such that, for all $\varepsilon \in (0, \varepsilon_\delta)$,*

$$\mathbb{P}_{\mathcal{J}_0}\{A(0, \tau_f(\varepsilon x)) - c\tau_f(\varepsilon x) > (1 - \delta)x\} \prod_{i \in \mathcal{J}_0} \mathbb{P}\{\mathbf{A}_i^r > \varepsilon x\}$$
$$\lesssim P_{\mathcal{J}_0}((1 - \delta)x) \prod_{i \in \mathcal{J}_1} p_i,$$

*with $P_{\mathcal{J}_0}((1 - \delta)x)$ as in (6.1).*

LEMMA 6.5 (Lower bound). *There exists an $\varepsilon > 0$ such that*

$$\mathbb{P}_{\mathcal{J}_0}\{A(0, \tau_f(\varepsilon x)) - c\tau_f(\varepsilon x) > x\} \prod_{i \in \mathcal{J}_0} \mathbb{P}\{\mathbf{A}_i^r > \varepsilon x\}$$
$$\lesssim P_{\mathcal{J}_0}(x) \prod_{i \in \mathcal{J}_1} p_i,$$

*with $P_{\mathcal{J}_0}(x)$ as in (6.1).*

The proofs of these lemmas are quite technical, and are deferred to Appendices A and B. A brief sketch of the proofs is given at the end of this section.

We now have gathered all the ingredients for the proof of Theorem 6.1.

PROOF OF THEOREM 6.1. The lower bound in Theorem 6.2 may be written as

$$\mathbb{P}\{A(0, \tau_f(\varepsilon x)) - c\tau_f(\varepsilon x) > x\}$$

$$= \sum_{\mathcal{J}_0 \subseteq \{1, \ldots, N\}} \mathbb{P}_{\mathcal{J}_0}\{A(0, \tau_f(\varepsilon x)) - c\tau_f(\varepsilon x) > x\} \mathbb{P}\{D_{\mathcal{J}_0}(\varepsilon x)\}.$$

Note that

$$\mathbb{P}\{D_{\mathcal{J}_0}(\varepsilon x)\} \sim \prod_{i \in \mathcal{J}_0} p_i \mathbb{P}\{\mathbf{A}_i^r > \varepsilon x\}.$$

Using Lemma 6.5, we then obtain

$$\mathbb{P}\{A(0, \tau_f(\varepsilon x)) - c\tau_f(\varepsilon x) > x\}$$

$$\gtrsim \left(\prod_{j=1}^N p_j\right) \sum_{\mathcal{J}_0 \subseteq \{1, \ldots, N\}} P_{\mathcal{J}_0}(x).$$

Similarly, using Lemma 6.4,

$$\mathbb{P}_{\mathcal{J}_0}\{A(0, \tau_f(\varepsilon x)) - c\tau_f(\varepsilon x) > (1 - \delta)x\}$$

$$\lesssim \left(\prod_{j=1}^N p_j\right) \sum_{\mathcal{J}_0 \subseteq \{1, \ldots, N\}} P_{\mathcal{J}_0}((1 - \delta)x).$$

Theorem 6.2 then gives

$$\left(\prod_{j=1}^N p_j\right) \sum_{\mathcal{J}_0 \subseteq \{1, \ldots, N\}} P_{\mathcal{J}_0}(x)$$

$$\lesssim \mathbb{P}\{\mathbf{V}^c > x\}$$

$$\lesssim \left(\prod_{j=1}^N p_j\right) \sum_{\mathcal{J}_0 \subseteq \{1, \ldots, N\}} P_{\mathcal{J}_0}((1 - \delta)x),$$

which implies Theorem 6.1, since $P_{\mathcal{J}_0}(x) \in \mathcal{R}$ as will be shown in Theorem 6.3. $\square$

In preparation for the proofs of Lemmas 6.4 and 6.5, we give a convenient representation for $A(0, \tau_f) - c\tau_f$ under the event $D_{\mathcal{J}_0}(\varepsilon x)$.

LEMMA 6.6. *Under the event* $D_{\mathcal{J}_0}(\varepsilon x)$, $A(0, \tau_f) - c\tau_f$ *can be represented as*

$$A(0, \tau_f) - c\tau_f$$

$$= \min\left\{\min_{i \in \mathcal{J}_0} \mathbf{F}_i, \min_{i \in \mathcal{J}_1} \mathbf{G}_i\right\},$$

*where $\mathcal{J}_1 = \mathcal{J} \setminus \mathcal{J}_0$. The random variables $\mathbf{F}_i$ and $\mathbf{G}_i$ are given by*

$$\mathbf{F}_i = (r - c)\bar{\mathbf{A}}_i^r(\varepsilon x)$$

$$- \sum_{k \in \mathcal{J}_1} r_k \left[ \mathbf{B}_k \mathbf{A}_k^r(\varepsilon x) + (1 - \mathbf{B}_k)[\mathbf{A}_k(\varepsilon x) + \mathbf{U}_k^r] + \sum_{j=1}^{N_k(\varepsilon x)} \mathbf{U}_{kj} \right],$$

$$\mathbf{G}_i = (r - c)\bar{\mathbf{A}}_i(\varepsilon x)$$

$$+ (r - c) \left[ \mathbf{B}_i \mathbf{A}_i^r(\varepsilon x) + (1 - \mathbf{B}_i)\mathbf{A}_i(\varepsilon x) + \sum_{j=1}^{N_i(\varepsilon x)} \mathbf{A}_{ij}(\varepsilon x) \right]$$

$$- d_i \left[ (1 - \mathbf{B}_i)\mathbf{U}_i^r + \sum_{j=1}^{N_i(\varepsilon x)} \mathbf{U}_{ij} \right] - \sum_{k \in \mathcal{J}_1 \setminus \{i\}} r_k \left[ (1 - \mathbf{B}_k)\mathbf{U}_k^r + \sum_{j=1}^{N_k(\varepsilon x)} \mathbf{U}_{kj} \right].$$

*Here $\bar{\mathbf{A}}_i(\varepsilon x) = \mathbf{A}_i | \mathbf{A}_i > \varepsilon x$, $\bar{\mathbf{A}}_i^r(\varepsilon x) = \mathbf{A}_i^r | \mathbf{A}_i^r > \varepsilon x$, $\mathbf{A}_{ij}(\varepsilon x) \stackrel{d}{=} \mathbf{A}_{ij} | \mathbf{A}_{ij} \leq \varepsilon x$ and $\mathbf{A}_i^r(\varepsilon x) \stackrel{d}{=} \mathbf{A}_i^r | \mathbf{A}_i^r \leq \varepsilon x$.*

PROOF. Under the event $D_{\mathcal{J}_0}(\varepsilon x)$, the random variables $\tau_{s,i}$, $i \in \mathcal{J}_1$, can be represented as

$$\tau_{s,i} = \mathbf{B}_i \mathbf{A}_i^r(\varepsilon x) + (1 - \mathbf{B}_i)[\mathbf{U}_i^r + \mathbf{A}_i(\varepsilon x)] + \sum_{j=1}^{N_i(\varepsilon x)} [\mathbf{U}_{ij} + \mathbf{A}_{ij}(\varepsilon x)], \qquad i \in \mathcal{J}_1.$$

Combined with the identities

$$A_i(0, \tau_{s,i}) = r_i \left[ \mathbf{B}_i \mathbf{A}_i^r(\varepsilon x) + (1 - \mathbf{B}_i)\mathbf{A}_i(\varepsilon x) + \sum_{j=1}^{N_i(\varepsilon x)} \mathbf{A}_{ij}(\varepsilon x) \right],$$

$$\tau_f = \min \left\{ \min_{i \in \mathcal{J}_0} \bar{\mathbf{A}}_i^r(\varepsilon x), \min_{i \in \mathcal{J}_1} \{ \bar{\mathbf{A}}_i(\varepsilon x) + \tau_{s,i} \} \right\},$$

$$A_i(\tau_{s,i}, \tau_f) = r_i(\tau_f - \tau_{s,i}),$$

the representation for $A(0, \tau_f) - c\tau_f$ then easily follows. □

We now give a brief sketch of the proofs of Lemmas 6.4 and 6.5. Both rely on the above representation for $A(0, \tau_f) - c\tau_f$ in terms of the variables $\mathbf{F}_i$ and $\mathbf{G}_i$. The proofs of the lemmas have a similar structure.

(i) The expressions for $\mathbf{F}_i$ and $\mathbf{G}_i$ are quite complicated, so an attempt to obtain the exact joint distribution does not seem promising. Therefore, the first step is to show that all random variables $\mathbf{A}_{ij}(\varepsilon x)$ and $\mathbf{U}_{ij}$ can be replaced by their means.

(ii) The above point indicates that $\mathbf{F}_i$ and $\mathbf{G}_i$ may be approximated as follows:

$$\mathbf{F}_i \approx (r-c)\bar{\mathbf{A}}_i^r(\varepsilon x) + \sum_{k \in \mathcal{J}_1} r_k \mathbb{E}\{\mathbf{U}_k\} N_k(\varepsilon x),$$

$$\mathbf{G}_i \approx (r-c)\bar{\mathbf{A}}_i(\varepsilon x) + [(r-c)\mathbb{E}\{\mathbf{A}_i\} - d_i \mathbb{E}\{\mathbf{U}_i\}] N_i(\varepsilon x)$$
$$- \sum_{k \in \mathcal{J}_1 \setminus \{i\}} r_k \mathbb{E}\{\mathbf{U}_k\} N_k(\varepsilon x).$$

It will be useful to keep these approximations in mind. The formulas in Appendixes A and B look much more cumbersome by the appearance of many additional but small constants.

(iii) The only random variables appearing in the above expressions are $\bar{\mathbf{A}}_i(\varepsilon x)$, $\mathbf{B}_i^r(\varepsilon x)$ and $N_i(\varepsilon x)$, of which the distributions are known. What thus remains is a straightforward computation.

The first point causes the most technical difficulties. It requires a separate treatment in the proofs of Lemmas 6.4 and 6.5. Details may be found in the Appendices.

6.4. *Asymptotic behavior of $P_{\mathcal{J}_0}(x)$ and $\mathbb{P}\{\mathbf{V}^c > x\}$.* In this section we give an asymptotic characterization of $P_{\mathcal{J}_0}(x)$, which may be useful for further analysis. In particular, we establish that $P_{\mathcal{J}_0}(x)$ and $\mathbb{P}\{\mathbf{V}^c > x\}$ are both regularly varying.

Assume that $\mathcal{J}_0$ is a proper subset of $\mathcal{J}$, observing

$$P_{\mathcal{J}}(x) = \prod_{i \in \mathcal{J}} \mathbb{P}\left\{\mathbf{A}_i^r > \frac{x}{r-c}\right\}.$$

For every set $\mathcal{J}_0$, define the $|\mathcal{J}_1|$-vector $g$ by

$$g := \left(\frac{r_j - \rho_j}{r-c}\right)_{j \in \mathcal{J}_1}.$$

Let $G$ be a (square) matrix with identical rows $g$, and let $\bar{G} := G - I$, with $I$ the identity matrix of dimension $|\mathcal{J}_1|$.

It can easily be shown that $\bar{G}$ is invertible; denote its inverse by $H$. A straightforward computation yields $H = \frac{1}{ge-1}G - I$, with $e = (1, \ldots, 1)$, which implies that $gH = \frac{1}{ge-1}g$. A further straightforward computation shows $|\bar{G}| = eg - 1$.

Define $y = (y_i)_{\mathcal{J}_1}$ and $dy = \prod_{i \in \mathcal{J}_1} dy_i$. Then we may write

$$P_{\mathcal{J}_0}(x) = \frac{1}{\prod_{i \in \mathcal{J}_1} \mathbb{E}\{\mathbf{A}_i\}}$$
$$\times \int_{y \geq 0} \prod_{i \in \mathcal{J}_1} \mathbb{P}\left\{\mathbf{A}_i > (\bar{G}y)_i + \frac{x}{r-c}\right\} \prod_{i \in \mathcal{J}_0} \mathbb{P}\left\{\mathbf{A}_i^r > gy + \frac{x}{r-c}\right\} dy.$$

If we integrate w.r.t. $z := \bar{G}y$ (note that $\bar{G}$ is a positive matrix), then we obtain [defining $\mathbf{A}_{\mathcal{J}_1} = (\mathbf{A}_i^r)_{i \in \mathcal{J}_1}$]

$$P_{\mathcal{J}_0}(x) = \frac{1}{|\bar{G}| \prod_{i \in \mathcal{J}_1} \mathbb{E}\{\mathbf{A}_i\}}$$

$$\times \int_{z \geq 0} \prod_{i \in \mathcal{J}_1} \mathbb{P}\left\{\mathbf{A}_i > z_i + \frac{x}{r-c}\right\} \prod_{i \in \mathcal{J}_0} \mathbb{P}\left\{\mathbf{A}_i^r > gHz + \frac{x}{r-c}\right\} dz$$

$$= \frac{1}{eg-1} \int_{z \geq 0} \prod_{i \in \mathcal{J}_0} \mathbb{P}\left\{\mathbf{A}_i^r > \frac{1}{eg-1} gz + \frac{x}{r-c}\right\} \prod_{i \in \mathcal{J}_1} d\mathbb{P}\left\{\mathbf{A}_i^r \leq z_i + \frac{x}{r-c}\right\}$$

$$= \frac{1}{eg-1} \mathbb{P}\left\{\mathbf{A}_i^r \geq \frac{x}{r-c}, i \in \mathcal{J};\right.$$

$$\left. \mathbf{A}_k^r - \frac{x}{r-c} \geq \frac{1}{eg-1} g\left(\mathbf{A}_{\mathcal{J}_1}^r - e\frac{x}{r-c}\right), k \in \mathcal{J}_1\right\}.$$

We conclude that $P_{\mathcal{J}_0}(x)$ can be written (up to a constant) as the probability that $(\mathbf{A}_i^r)_{i \in \mathcal{J}}$ belongs to a certain set. We now show that $P_{\mathcal{J}_0}(x)$ is regularly varying of index $-\mu$ [recall that $\mu = \sum_{i=1}^{N}(\nu_i - 1)$]. If $\mathbf{A}_i$ is regularly varying of index $-\nu_i < -1$, then it is well known and elementary to show that

$$\mathbb{P}\left\{\frac{\mathbf{A}_i^r - \gamma x}{x} > y \,|\, \mathbf{A}_i^r > \gamma x\right\} \to \left(1 + \frac{y}{\gamma}\right)^{1-\nu_i},$$

as $x \to \infty$. Let $\mathbf{Z}_i$ be a random variable with the above limiting distribution, with $\gamma = \frac{1}{r-c}$ such that the $\mathbf{Z}_i$, $i \in \mathcal{J}_1$, are independent. The above computations are summarized in the following theorem.

THEOREM 6.3.

$$P_{\mathcal{J}_0}(x) \sim \kappa_{\mathcal{J}_0} \prod_{i=1}^{N} \mathbb{P}\left\{\mathbf{A}_i^r > \frac{x}{r-c}\right\},$$

with $\kappa_{\mathcal{J}} = 1$ and

$$\kappa_{\mathcal{J}_0} = \frac{1}{eg-1} \mathbb{P}\left\{\mathbf{Z}_i \geq \frac{1}{eg-1} g\mathbf{Z}_{\mathcal{J}_1}, i \in \mathcal{J}_0\right\}$$

if $\mathcal{J}_0$ is a proper subset of $\mathcal{J}$. In particular, $P_{\mathcal{J}_0}(x)$ is regularly varying of index $-\mu$.

Combining Theorems 6.1 and 6.3, we obtain:

THEOREM 6.4.

$$\mathbb{P}\{\mathbf{V}^c > x\} \sim \kappa \prod_{i=1}^{N} p_i \mathbb{P}\left\{\mathbf{A}_i^r > \frac{x}{r-c}\right\},$$

*with*

$$\kappa = \sum_{\mathcal{J}_0 \subseteq \{1, \ldots, N\}} \kappa_{\mathcal{J}_0}.$$

*In particular, $\mathbb{P}\{\mathbf{V}^c > x\}$ is regularly varying of index $-\mu$.*

The above theorem is used in proving the reduced-load equivalence (see Section 5), and may be potentially useful for computational purposes.

In particular, in the case of two On–Off flows, the computation of $\kappa$ is as difficult as the computation of $\kappa_1$ and $\kappa_2$. Using the probabilistic interpretation of these constants readily leads to an integral expression, which can be solved explicitly when both $\nu_1$ and $\nu_2$ are integer-valued. We omit the details.

**7. $K$ heterogeneous classes: proofs.** In this section we provide the proofs of the results in Section 3.5 for the case with $K$ heterogeneous classes of On–Off flows. In particular, we present a proof of Theorem 3.3.

We start with the regime where we first let $x \to \infty$ and then $n \to \infty$. For every $n$ we have, using Theorem 3.2,

$$\lim_{x \to \infty} \frac{\mathbb{P}\{\mathbf{V}^{(n)} > nx\}}{\log x} = -\mu^{(n)},$$

with $\mu^{(n)}$ denoting the optimal value of the criterion function of the associated knapsack problem. It thus remains to be shown that

$$(7.1) \qquad \lim_{n \to \infty} \frac{\mu^{(n)}}{n\mu} = 1.$$

First observe that the optimal value of the continuous relaxation of the knapsack problem is $n\mu$, yielding a lower bound for $\mu^{(n)}$. On the other hand, the continuous relaxation may be used to construct a feasible solution of the knapsack problem. Take (use the notation of Section 3.5) $q_k = n_k = na_k$ for $k < l$, $q_k = n_k = 0$ for $k > l$, and $q_l = |n_l| + 1$. This is a feasible solution with a value at most $n\mu + \gamma_l$, giving an upper bound for $\mu^{(n)}$. In conclusion, we have

$$n\mu \le \mu^{(n)} \le n\mu + \gamma_l,$$

from which (7.1) directly follows.

We now turn to the regime where we first let $n \to \infty$ and then $x \to \infty$ (i.e., the many-sources regime). Define the "decay rate"

$$I(x) := -\lim_{n \to \infty} \frac{1}{n} \log \mathbb{P}\{\mathbf{V}^{(n)} > nx\}.$$

It needs to be shown that $I(x) \sim \mu \log x$ as $x \to \infty$.

The above decay rate equals ([8], page 300)

$$I(x) = \inf_{t \geq 0} \sup_{\theta} \left( \theta(x+t) - \sum_{k=1}^{K} a_k \log \mathbb{E}\{e^{\theta A_k(t)}\} \right),$$

with $A_k(t) := A_k(0, t)$ representing the amount of traffic generated by a single class-$k$ flow in a time interval of length $t$ in steady state. Replacing $\theta$ by $\theta(\log t)/t$, we obtain an alternative variational problem,

$$\inf_{t \geq 0} \log t \cdot J_t\left(\frac{x}{t} + 1\right)$$

(7.2)

$$\text{where } J_t(x) := \sup_{\theta} \left( \theta x - \sum_{k=1}^{K} a_k \frac{\log \mathbb{E}\{e^{\theta(\log t)A_k(t)/t}\}}{\log t} \right),$$

for $x \in (0, \sum_{k=1}^{K} a_k r_k)$. The latter variational problem allows direct asymptotic analysis ($x \to \infty$) as in [29], which yields Theorem 7.1.

First, however, we state an auxiliary lemma. Recall that $\sigma_k = \sum_{m=1}^{k-1} a_m r_m + \sum_{m=k}^{K} a_m \rho_m$, and that the $K$ classes are indexed in nondecreasing order of the ratios $\gamma_k = (\nu_k - 1)/(r_k - \rho_k)$.

LEMMA 7.1. *For $\theta \geq 0$,*

$$\lim_{t \to \infty} \frac{\log \mathbb{E}\{e^{\theta(\log t)A_k(t)/t}\}}{\log t} = \max\{\theta \rho_k, \theta r_k - \nu_k + 1\},$$

*so that the cumulant function of the superposition is piecewise linear:*

$$\sum_{k=1}^{K} a_k \lim_{t \to \infty} \frac{\log \mathbb{E}\{e^{\theta(\log t)A_k(t)/t}\}}{\log t} = \sum_{k=1}^{K} a_k \max\{\theta \rho_k, \theta r_k - \nu_k + 1\}.$$

*Further,*

$$(7.3) \quad \lim_{t \to \infty} J_t(x) = \gamma_{l(x)} x - \sum_{k=1}^{l(x)-1} a_k(\gamma_{l(x)} r_k - \nu_k + 1) - \sum_{k=l(x)}^{K} a_k \gamma_{l(x)} \rho_k,$$

*for $x \in (0, \sum_{k=1}^{K} a_k r_k)$, where $l(x)$ is such that $x \in (\sigma_{l(x)-1}, \sigma_{l(x)})$.*
*The function $\lim_{t \to \infty} J_t(\cdot)$ is increasing.*

The proof of the above lemma is analogous to that of Theorem 3.6 and Lemma 3.7 of [29].

THEOREM 7.1 (Large-buffer asymptotics).

$$\lim_{x \to \infty} \frac{I(x)}{\log x} = \mu,$$

*with $\mu = \sum_{k=1}^{l-1} a_k(\nu_k - 1) + (1 - \sigma_{l-1})\gamma_l$ and $l := l(1)$.*

Proof.   The proof consists of deriving an upper bound and a lower bound which asymptotically coincide.

(Upper bound.)   Using the representation (7.2),

$$\limsup_{x \to \infty} \frac{I(x)}{\log x} = \limsup_{x \to \infty} \inf_{t > 0} \frac{\log t}{\log x} J_t \left( \frac{x}{t} + 1 \right).$$

Substituting $t = x/s$, $s \in (0, \sum_{k=1}^{K} a_k r_k - 1)$, to obtain an upper bound, and using (7.3),

$$\limsup_{x \to \infty} \inf_{t > 0} \frac{\log t}{\log x} J_t \left( \frac{x}{t} + 1 \right)$$

$$\leq \limsup_{x \to \infty} \frac{\log(x/s)}{\log x} J_{x/s}(s + 1)$$

$$\leq \limsup_{x \to \infty} \frac{\log(x/s)}{\log x} \limsup_{x \to \infty} J_{x/s}(s + 1)$$

$$\leq \limsup_{x \to \infty} J_{x/s}(s + 1)$$

$$= \gamma_{l(s+1)}(s + 1)$$

$$\qquad - \sum_{k=1}^{l(s+1)-1} (a_k \gamma_{l(s+1)} r_k - \nu_k + 1) - \sum_{k=l(s+1)}^{K} a_k \gamma_{l(s+1)} \rho_k.$$

The above inequality holds for any $s \in (0, \sum_{k=1}^{K} a_k r_k - 1)$. According to Lemma 7.1, the last term is increasing in $s + 1$. Letting $s \downarrow 0$ to obtain the sharpest possible upper bound, we obtain

$$\limsup_{x \to \infty} \frac{I(x)}{\log x} \leq \gamma_l - \sum_{k=1}^{l-1} a_k (\gamma_l r_k - \nu_k + 1) - \sum_{k=l}^{K} a_k \gamma_l \rho_k = \mu.$$

(Lower bound.)   Using the representation (7.2), and taking $\theta = \gamma_l$, we obtain the lower bound

$$I(x) = \inf_{t \geq 0} \log t \cdot \sup_{\theta} \left( \theta \left( \frac{x}{t} + 1 \right) - \sum_{k=1}^{K} a_k \frac{\log \mathbb{E}\{e^{\theta(\log t) A_k(t)/t}\}}{\log t} \right)$$

$$\geq \inf_{t \geq 0} \log t \cdot \left( \gamma_l \left( \frac{x}{t} + 1 \right) - \sum_{k=1}^{K} a_k \frac{\log \mathbb{E}\{e^{\gamma_l(\log t) A_k(t)/t}\}}{\log t} \right).$$

The optimizing value of $t$ in the above variational problem is *at least* linear in $x$, for large $x$. Formally, there exists a $d$ such that the above infimum need be taken only over $t > dx$, for large $x$. This may be proven analogously to case (iii) of [14], page 258. Thus,

$$I(x) \geq \inf_{t > dx} \log t \cdot \left( \gamma_l \left( \frac{x}{t} + 1 \right) - \sum_{k=1}^{K} a_k \frac{\log \mathbb{E}\{e^{\gamma_l(\log t) A_k(t)/t}\}}{\log t} \right).$$

Using (7.3), we find that for any $\varepsilon > 0$, and $x$ large enough, we have, for all $t > dx$,

$$\sum_{k=1}^{K} a_k \frac{\log \mathbb{E}\{e^{\gamma_l \log t A_k(t)/t}\}}{\log t} \le (1+\varepsilon) \sum_{k=1}^{K} a_k \max\{\gamma_l \rho_k, \gamma_l r_k - \nu_k + 1\}.$$

Thus,

$$\liminf_{x \to \infty} \frac{I(x)}{\log x}$$

$$\ge \liminf_{x \to \infty} \inf_{t > dx} \frac{\log t}{\log x}\left(\gamma_l\left(\frac{x}{t}+1\right)\right.$$

$$\left. - (1+\varepsilon)\sum_{k=1}^{K} a_k \max\{\gamma_l \rho_k, \gamma_l r_k - \nu_k + 1\}\right)$$

$$\ge \liminf_{x \to \infty} \inf_{t > dx} \frac{\log t}{\log x} \inf_{t > dx}\left(\gamma_l\left(\frac{x}{t}+1\right)\right.$$

$$\left. - (1+\varepsilon)\sum_{k=1}^{K} a_k \max\{\gamma_l \rho_k, \gamma_l r_k - \nu_k + 1\}\right)$$

$$\ge \gamma_l - (1+\varepsilon)\sum_{k=1}^{K} a_k \max\{\gamma_l \rho_k, \gamma_l r_k - \nu_k + 1\}.$$

Letting $\varepsilon \downarrow 0$, we obtain

$$\liminf_{x \to \infty} \frac{I(x)}{\log x}$$

$$\ge \gamma_l - \sum_{k=1}^{K} a_k \max\{\gamma_l \rho_k, \gamma_l r_k - \nu_k + 1\}$$

$$= \gamma_l - \sum_{k=1}^{K} a_k(\gamma_l \rho_k + \max\{0, \gamma_l(r_k - \rho_k) - \nu_k + 1\})$$

$$= \gamma_l - \sum_{k=1}^{K} a_k(\gamma_l \rho_k + \max\{0, (\gamma_l - \gamma_k)(r_k - \rho_k)\})$$

$$= \gamma_l - \sum_{k=1}^{l-1} a_k(\gamma_l r_k - \nu_k + 1) - \sum_{k=l}^{K} a_k \gamma_l \rho_k$$

$$= \mu. \qquad \qquad \qquad \square$$

As shown above, Theorem 3.3 implies that the limits $x \to \infty$ and $n \to \infty$ commute, as long as one considers "rough" (i.e., logarithmic) asymptotics. However, in case of "more refined" asymptotics, the limits do not necessarily commute. This may be seen as follows. Consider the case of $n$ homogeneous On–Off flows with Pareto($\nu$) distributed On periods. In [28], it is proven that

$$\lim_{x \to \infty} \lim_{n \to \infty} \frac{1}{n} \log \mathbb{P}\{\mathbf{V}^{(n)} > nx\} + (\nu - 1)\left(\frac{c - \rho}{r - \rho}\right) \log(x \log x) = H,$$

for some constant $H \in (0, \infty)$. Now reverse the limits. Denote by $k_n$ the number of flows sending at peak rate in the reduced-load approximation (in the notation of Section 3.4, we have $k_n = N^*$):

$$k_n := \left\lceil \frac{nc - n\rho}{r - \rho} \right\rceil.$$

Now with Theorem 3.1, we have, for any finite $n$ and $x \to \infty$,

$$\mathbb{P}\{\mathbf{V}^{(n)} > nx\} \sim f(n) x^{-(\nu - 1)k_n},$$

for some function $f(\cdot)$. Hence,

$$\lim_{x \to \infty} \frac{1}{n} \log \mathbb{P}\{\mathbf{V}^{(n)} > nx\} + (\nu - 1)\left(\frac{c - \rho}{r - \rho}\right) \log(x \log x)$$

$$= \log f(n) + \lim_{x \to \infty} (\nu - 1)\left(\frac{k_n}{n} - \frac{c - \rho}{r - \rho}\right) \log x - (\nu - 1)\frac{c - \rho}{r - \rho} \log \log x.$$

Since this limit does not exist in $\mathbb{R}$, we conclude that the limits do not necessarily commute.

## 8. Concluding remarks.

We have characterized the asymptotic behavior of the workload distribution in a fluid queue fed by multiple heavy-tailed On–Off flows. The results extend previous work, like the bounds derived in [15], and the exact asymptotics in [9, 22] which rely on strong peak-rate conditions. As a by-product, the proofs lead to several important insights like the extension of the reduced-load equivalence established in [1] (see Section 5), and a detailed understanding of the typical overflow behavior (see Section 6). In the analysis, we excluded the case where the drift may be zero during the path to overflow (see Section 3.1 for a brief discussion), which appears particularly interesting from a theoretical perspective.

There are several other interesting topics for further research. We expect that the methodology of Section 6 is also suitable to study other similar problems, such as fluid queues with $M/G/\infty$ input, multiserver queues, and Generalized Processor Sharing queues. A further avenue for research is the extension of the results to the case of On–Off flows with more general subexponential On periods, for example, Weibull. Partial results in [1] indicate that the typical overflow behavior may then actually be quite different.

## APPENDIX A

**Proof of Lemma 6.4.**

LEMMA 6.4 (Upper bound).   *For any $\delta > 0$, there exists an $\varepsilon_\delta > 0$ such that, for all $\varepsilon \in (0, \varepsilon_\delta)$,*

$$\mathbb{P}_{\mathcal{J}_0}\{A(0, \tau_f(\varepsilon x)) - c\tau_f(\varepsilon x) > (1 - \delta)x\} \prod_{i \in \mathcal{J}_0} \mathbb{P}\{\mathbf{A}_i^r > \varepsilon x\}$$

$$\lesssim P_{\mathcal{J}_0}((1 - \delta)x) \prod_{i \in \mathcal{J}_1} p_i,$$

*with $P_{\mathcal{J}_0}((1 - \delta)x)$ as in (6.1).*

PROOF.   As mentioned earlier, the first step is to replace all random variables $\mathbf{A}_{ij}$ and $\mathbf{U}_{ij}$ by their means. Let $\bar{\delta}$ and $\tilde{\delta}$ be two $|\mathcal{J}_1|$-vectors, of which the elements are positive but arbitrarily small. Note that, for fixed $\mathcal{J}_0$,

$$\mathbf{F}_i \leq (r - c)\bar{\mathbf{A}}_i^r(\varepsilon x) - \sum_{k \in \mathcal{J}_1} r_k N_k(\varepsilon x)[\mathbb{E}\{\mathbf{U}_k\} - \bar{\delta}_k]$$

$$+ \sum_{k \in \mathcal{J}_1} r_k \sum_{j=1}^{N_k(\varepsilon x)} [\mathbb{E}\{\mathbf{U}_k\} - \bar{\delta}_k - \mathbf{U}_{kj}],$$

$$\mathbf{G}_i \leq (r - c)\bar{\mathbf{A}}_i(\varepsilon x) + (r - c)\varepsilon x + (r - c)N_i(\varepsilon x)[\mathbb{E}\{\mathbf{A}_i\} + \tilde{\delta}_i]$$

$$+ (r - c)\sum_{j=1}^{N_i(\varepsilon x)}[\mathbf{A}_{ij}(\varepsilon x) - \mathbb{E}\{\mathbf{A}_i\} - \tilde{\delta}_i] - d_i N_i(\varepsilon x)[\mathbb{E}\{\mathbf{U}_i\} - \tilde{\delta}_i]$$

$$+ d_i \sum_{j=1}^{N_i(\varepsilon x)} [\mathbb{E}\{\mathbf{U}_i\} - \tilde{\delta}_i - \mathbf{U}_{ij}] \sum_{k \in \mathcal{J}_1 \setminus \{i\}} r_k N_k(\varepsilon x)[\mathbb{E}\{\mathbf{U}_k\} - \bar{\delta}_k]$$

$$+ \sum_{k \in \mathcal{J}_1 \setminus \{i\}} r_k \sum_{j=1}^{N_k(\varepsilon x)} [\mathbb{E}\{\mathbf{U}_k\} - \bar{\delta}_k - \mathbf{U}_{kj}].$$

Define the event $E_1(\gamma, \bar{\delta}, \tilde{\delta}, \varepsilon, x)$ by

$$\left\{ \sum_{j=1}^{N_i(\varepsilon x)} [\mathbb{E}\{\mathbf{U}_i\} - \min\{\bar{\delta}_i, \tilde{\delta}_i\} - \mathbf{U}_{ij}] \leq \gamma x/(2r), i \in \mathcal{J}_1 \right\}$$

$$\cup \left\{ \sum_{j=1}^{N_i(\varepsilon x)} [\mathbf{A}_{ij}(\varepsilon x) - \mathbb{E}\{\mathbf{A}_i\} - \min\{\bar{\delta}_i, \tilde{\delta}_i\}] \right.$$

$$\leq \gamma x/(2r) - (r-c)\varepsilon x, i \in \mathcal{J}_1\Big\}.$$

A straightforward application of Lemma 6.2 (analogously to the proofs of Propositions 6.2, 6.4 and 6.5) shows that, for any $\gamma, \bar{\delta}, \tilde{\delta} > 0$, there exists an $\varepsilon^* > 0$ such that, for all $\varepsilon \in (0, \varepsilon^*]$,

$$\text{(A.1)} \qquad \mathbb{P}_{\mathcal{J}_0}\{E_1(\gamma, \bar{\delta}, \tilde{\delta}, \varepsilon, x)^c\} = o(P(x)),$$

as $x \to \infty$ with $P(x) = \prod_{j=1}^N \mathbb{P}\{\mathbf{A}_j^r > x\}$, as defined earlier.

From (A.1) and Lemma 6.6, we conclude that, using the upper bounds for $\mathbf{F}_i$ and $\mathbf{G}_i$,

$$\mathbb{P}_{\mathcal{J}_0}\{A(0, \tau_f) - c\tau_f > (1 - \delta)x\}$$

$$= \mathbb{P}_{\mathcal{J}_0}\{A(0, \tau_f) - c\tau_f > (1 - \delta)x; E_1(\gamma, \bar{\delta}, \tilde{\delta}, \varepsilon, x)^c\}$$

$$\quad + \mathbb{P}_{\mathcal{J}_0}\{A(0, \tau_f) - c\tau_f > (1 - \delta)x; E_1(\gamma, \bar{\delta}, \tilde{\delta}, \varepsilon, x)\}$$

$$\leq \mathbb{P}\Bigg\{(r-c)\bar{\mathbf{A}}_i^r(\varepsilon x)$$

$$\quad - \sum_{k \in \mathcal{J}_1} r_k N_k(\varepsilon x)[\mathbb{E}\{\mathbf{U}_k\} - \bar{\delta}_k] > (1 - \gamma - \delta)x, i \in \mathcal{J}_0;$$

$$\quad (r-c)\bar{\mathbf{A}}_i(\varepsilon x) + (r-c)N_i(\varepsilon x)[\mathbb{E}\{\mathbf{A}_i\} + \tilde{\delta}_i]$$

$$\quad - d_i N_i(\varepsilon x)[\mathbb{E}\{\mathbf{U}_i\} - \tilde{\delta}_i]$$

$$\quad - \sum_{k \in \mathcal{J}_1 \setminus \{i\}} r_k N_k(\varepsilon x)[\mathbb{E}\{\mathbf{U}_k\} - \bar{\delta}_k] > (1 - \gamma - \delta)x, i \in \mathcal{J}_1\Bigg\}$$

$$\quad + o(P(x)).$$

The last probability equals [condition on $N_i(\varepsilon x)$, $i \in \mathcal{J}_1$]

$$\sum_{n_i \geq 1, i \in \mathcal{J}_1} \Bigg(\prod_{i \in \mathcal{J}_1} \mathbb{P}\{N_i(\varepsilon x) = n_i\}\Bigg)$$

$$\quad \times \mathbb{P}\Bigg\{(r-c)\bar{\mathbf{A}}_i^r(\varepsilon x)$$

$$\quad - \sum_{k \in \mathcal{J}_1} r_k[\mathbb{E}\{\mathbf{U}_k\} - \bar{\delta}_k]n_k > (1 - \gamma - \delta)x, i \in \mathcal{J}_0;$$

$$\quad (r-c)\bar{\mathbf{A}}_i(\varepsilon x) + (r-c)[\mathbb{E}\{\mathbf{A}_i\} + \tilde{\delta}_i]n_i$$

$$\quad - d_i[\mathbb{E}\{\mathbf{U}_i\} - \tilde{\delta}_i]n_i$$

$$- \sum_{k \in \mathcal{J}_1 \setminus \{i\}} r_k [\mathbb{E}\{\mathbf{U}_k\} - \bar{\delta}_k] n_k > (1 - \gamma - \delta)x, i \in \mathcal{J}_1 \bigg\}.$$

Deconditioning upon $\bar{\mathbf{A}}_i$ and $\bar{\mathbf{A}}_i^r$ (i.e., dividing by $\prod_{i \in \mathcal{J}_0} \mathbb{P}\{\mathbf{A}_i^r > \varepsilon x\} \prod_{i \in \mathcal{J}_1} \mathbb{P}\{\mathbf{A}_i > \varepsilon x\}$), and noting that $\mathbb{P}\{N_i(\varepsilon x) = n_i\} \leq \mathbb{P}\{\mathbf{A}_i > \varepsilon x\}$, we obtain that

$$\mathbb{P}_{\mathcal{J}_0}\{A(0, \tau_f) - c\tau_f > (1 - \delta)x\}$$

$$\times \prod_{i \in \mathcal{J}_0} \mathbb{P}\{\mathbf{A}_i^r > \varepsilon x\}$$

is upper bounded by [up to $o(P(x))$]

$$\sum_{n_i \geq 0, i \in \mathcal{J}_1} \bigg( \prod_{i \in \mathcal{J}_0} \mathbb{P}\bigg\{ (r - c)\mathbf{A}_i^r > (1 - \gamma - \delta)x$$

$$+ \sum_{k \in \mathcal{J}_1} r_k [\mathbb{E}\{\mathbf{U}_k\} - \bar{\delta}_k] n_k \bigg\} \bigg)$$

$$\times \prod_{i \in \mathcal{J}_1} \mathbb{P}\bigg\{ (r - c)\mathbf{A}_i > (1 - \gamma - \delta)x$$

$$+ [d_i \mathbb{E}\{\mathbf{U}_i\} - (r - c)\mathbb{E}\{\mathbf{A}_i\} - r_i \tilde{\delta}_i] n_i$$

$$+ \sum_{k \in \mathcal{J}_1 \setminus \{i\}} r_k [\mathbb{E}\{\mathbf{U}_k\} - \bar{\delta}_k] n_k \bigg\}.$$

It is important to note that this expression is independent of $\varepsilon$.

Since all probabilities appearing in the right-hand side are decreasing functions of $n_i$ (for $\bar{\delta}$ and $\tilde{\delta}$ small enough), the latter term is bounded by [with $y := (y_i)_{i \in \mathcal{J}_1}$ and $dy := \prod_{i \in \mathcal{J}_1} dy_i$]

$$\int_{y \geq 0} \mathbb{P}\bigg\{ (r - c)\mathbf{A}_i^r > (1 - \gamma - \delta)x + \sum_{k \in \mathcal{J}_1} r_k [\mathbb{E}\{\mathbf{U}_k\} - \bar{\delta}_k] y_k \bigg\}$$

$$\times \prod_{i \in \mathcal{J}_1} \mathbb{P}\bigg\{ (r - c)\mathbf{A}_i > (1 - \gamma - \delta)x$$

(A.2)
$$+ [d_i \mathbb{E}\{\mathbf{U}_i\} - (r - c)\mathbb{E}\{\mathbf{A}_i\} - r_i \tilde{\delta}_i] y_i$$

$$+ \sum_{k \in \mathcal{J}_1 \setminus \{i\}} r_k [\mathbb{E}\{\mathbf{U}_k\} - \bar{\delta}_k] y_k \bigg\} dy.$$

We will rewrite this expression in terms of $P_{\mathcal{J}_0}(x)$. Apply the change of variables $z_i := y_i / (\mathbb{E}\{\mathbf{A}_i\} + \mathbb{E}\{\mathbf{U}_i\})$. Redefine $\bar{\delta}_i := \bar{\delta}_i / (\mathbb{E}\{\mathbf{A}_i\} + \mathbb{E}\{\mathbf{U}_i\})$ and similarly $\tilde{\bar{\delta}}_i := \bar{\delta}_i / (\mathbb{E}\{\mathbf{A}_i\} + \mathbb{E}\{\mathbf{U}_i\})$. Note that

$$\frac{1}{\mathbb{E}\{\mathbf{A}_i\} + \mathbb{E}\{\mathbf{U}_i\}} = \frac{p_i}{\mathbb{E}\{\mathbf{A}_i\}} \quad \text{and} \quad r_i \frac{\mathbb{E}\{\mathbf{U}_i\}}{\mathbb{E}\{\mathbf{A}_i\} + \mathbb{E}\{\mathbf{U}_i\}} = r_i(1 - p_i) = r_i - \rho_i.$$

Then we obtain that (A.2) equals

$$\left( \prod_{i \in \mathcal{J}_1} \frac{p_i}{\mathbb{E}\{\mathbf{A}_i\}} \right)$$

$$\times \int_{z \geq 0} \prod_{i \in \mathcal{J}_0} \mathbb{P}\left\{ (r - c)\mathbf{A}_i^r > (1 - \gamma - \delta)x + \sum_{k \in \mathcal{J}_1} (r_k - \rho_k - \bar{\delta}_k)z_k \right\}$$

$$\times \prod_{i \in \mathcal{J}_1} \mathbb{P}\left\{ (r - c)\mathbf{A}_i > (1 - \gamma - \delta)x + (d_i - \rho_i - \tilde{\bar{\delta}}_i)z_i \right.$$

$$\left. + \sum_{k \in \mathcal{J}_1 \setminus \{i\}} (r_k - \rho_k - \bar{\delta}_k)z_k \right\} dz.$$

If we take $\tilde{\bar{\delta}}_i = \frac{d_i - \rho_i}{r_i - \rho_i} \bar{\delta}_i$ and integrate w.r.t. $z_i \frac{r_i - \rho_i - \bar{\delta}_i}{r_i - \rho_i}$, then we obtain

$$\left( \prod_{i \in \mathcal{J}_1} \frac{r_i - \rho_i}{r_i - \rho_i - \bar{\delta}_i} \frac{p_i}{\mathbb{E}\{\mathbf{A}_i\}} \right)$$

$$\times \int_{z \geq 0} \prod_{i \in \mathcal{J}_0} \mathbb{P}\left\{ (r - c)\mathbf{A}_i^r > (1 - \gamma - \delta)x + \sum_{k \in \mathcal{J}_1} (r_k - \rho_k)z_k \right\}$$

$$\times \prod_{i \in \mathcal{J}_1} \mathbb{P}\left\{ (r - c)\mathbf{A}_i > (1 - \gamma - \delta)x \right.$$

$$\left. + (d_i - \rho_i)z_i + \sum_{k \in \mathcal{J}_1 \setminus \{i\}} (r_k - \rho_k)z_k \right\} dz$$

$$= \prod_{i \in \mathcal{J}_1} p_i \frac{r_i - \rho_i}{r_i - \rho_i - \bar{\delta}_i} P_{\mathcal{J}_0}((1 - \gamma - \delta)x).$$

Together with the fact that $P_{\mathcal{J}_0}(\cdot)$ is regularly varying (see below), this completes the proof of the upper bound after dividing by $P_{\mathcal{J}_0}(x)$, letting $x \to \infty$, and noting that $\delta$, $\bar{\delta}$ and $\gamma$ may be chosen arbitrarily small. $\quad\square$

## APPENDIX B

LEMMA 6.5 (Lower bound).    *There exists an $\varepsilon > 0$ such that*

$$\mathbb{P}_{\mathcal{J}_0}\{A(0, \tau_f(\varepsilon x)) - c\tau_f(\varepsilon x) > x\} \prod_{i \in \mathcal{J}_0} \mathbb{P}\{\mathbf{A}_i^r > \varepsilon x\} \gtrsim P_{\mathcal{J}_0}(x) \prod_{i \in \mathcal{J}_1} p_i,$$

with $P_{\mathcal{J}_0}(x)$ as in (6.1).

PROOF.   As in Appendix A, the first step is to replace the random variables $\mathbf{A}_i(\varepsilon x)$ and $\mathbf{U}_i$ by their means. Adding and subtracting appropriate means, it is easy to see that, for fixed $\mathcal{J}_0$,

$$\mathbf{F}_i = (r - c)\bar{\mathbf{A}}_i^r(\varepsilon x) - \sum_{k \in \mathcal{J}_1} r_k N_k(\varepsilon x)[\mathbb{E}\{\mathbf{U}_k\} + \bar{\delta}_k]$$

$$+ \sum_{k \in \mathcal{J}_1} r_k \sum_{j=1}^{N_k(\varepsilon x)} [\mathbb{E}\{\mathbf{U}_k\} - \mathbf{U}_{kj} + \bar{\bar{\delta}}_k]$$

$$- \sum_{k \in \mathcal{J}_1} r_k[\mathbf{B}_k \mathbf{A}_k^r(\varepsilon x) + (1 - \mathbf{B}_k)(\mathbf{A}_k(\varepsilon x) + \mathbf{U}_k^r)],$$

$$\mathbf{G}_i = (r - c)\bar{\mathbf{A}}_i(\varepsilon x) + (r - c)[\mathbf{B}_i \mathbf{A}_i^r(\varepsilon x) + (1 - \mathbf{B}_i)\mathbf{A}_i(\varepsilon x)]$$

$$- d_i(1 - \mathbf{B}_i)\mathbf{U}_i^r - \sum_{k \in \mathcal{J}_1 \setminus \{i\}} r_k(1 - \mathbf{B}_k)\mathbf{U}_k^r$$

$$+ (r - c) \sum_{j=1}^{N_i(\varepsilon x)} [\mathbf{A}_{ij}(\varepsilon x) - \mathbb{E}\{\mathbf{A}_i\} + \tilde{\delta}_i]$$

$$+ (r - c)N_i(\varepsilon x)[\mathbb{E}\{\mathbf{A}_i\} - \tilde{\delta}_i]$$

$$- d_i N_i(\varepsilon x)[\mathbb{E}\{\mathbf{U}_i\} + \tilde{\bar{\delta}}_i]$$

$$+ d_i \sum_{j=1}^{N_i(\varepsilon x)} [\mathbb{E}\{\mathbf{U}_i\} - \mathbf{U}_{ij} + \tilde{\delta}_i]$$

$$- \sum_{k \in \mathcal{J}_1 \setminus \{i\}} r_k N_k(\varepsilon x)[\mathbb{E}\{\mathbf{U}_k\} + \bar{\delta}_k]$$

$$+ \sum_{k \in \mathcal{J}_1 \setminus \{i\}} r_k \sum_{j=1}^{N_k(\varepsilon x)} [\mathbb{E}\{\mathbf{U}_k\} - \mathbf{U}_{kj} + \bar{\delta}_k].$$

Define the event $E_2(\gamma, \bar{\delta}, \tilde{\delta}, \varepsilon, x)$ by

$$\left\{ \sum_{j=1}^{N_i(\varepsilon x)} [\mathbb{E}\{\mathbf{U}_i\} - \mathbf{U}_{ij} + \min\{\bar{\delta}_i, \tilde{\delta}_i\}] \geq -\gamma x/(3r), i \in \mathcal{J}_1 \right\}$$

$$\cup \left\{ \sum_{j=1}^{N_i(\varepsilon x)} [\mathbf{A}_{ij}(\varepsilon x) - \mathbb{E}\{\mathbf{A}_i\} + \min\{\bar{\delta}_i, \tilde{\delta}_i\}] \geq -\gamma x/(3r), i \in \mathcal{J}_1 \right\}$$

$$\cup \left\{ \sum_{k \in \mathcal{J}_1} [\mathbf{B}_k \mathbf{A}_k^r(\varepsilon x) + (1 - \mathbf{B}_k)(\mathbf{A}_k(\varepsilon x) + \mathbf{U}_k^r)] \leq \gamma x/(3r) \right\}.$$

We have the lower bound

$$\mathbb{P}_{\mathcal{J}_0}\{A(0, \tau_f) - c\tau_f > x\}$$

$$= \mathbb{P}_{\mathcal{J}_0}\{\mathbf{F}_i > x, i \in \mathcal{J}_0; \mathbf{G}_i > x, i \in \mathcal{J}_1\}$$

$$\geq \mathbb{P}\{\mathbf{F}_i > x, i \in \mathcal{J}_0; \mathbf{G}_i > x, i \in \mathcal{J}_1; E_2(\gamma, \bar{\delta}, \tilde{\delta}, \varepsilon, x)\}$$

$$\geq \mathbb{P}\bigg\{ (r - c)\bar{\mathbf{A}}_i^r(\varepsilon x) - \sum_{k \in \mathcal{J}_1} r_k N_k(\varepsilon x)[\mathbb{E}\{\mathbf{U}_k\} + \bar{\delta}_k] > (1 + \gamma)x, i \in \mathcal{J}_0;$$

$$(r - c)\bar{\mathbf{A}}_i(\varepsilon x) + (r - c)N_i(\varepsilon x)[\mathbb{E}\{\mathbf{A}_i\} - \tilde{\delta}_i] - d_i N_i(\varepsilon x)[\mathbb{E}\{\mathbf{U}_i\} + \tilde{\delta}_i]$$

$$- \sum_{k \in \mathcal{J}_1 \setminus \{i\}} r_k N_k(\varepsilon x)[\mathbb{E}\{\mathbf{U}_k\} + \bar{\delta}_k] > (1 + \gamma)x, i \in \mathcal{J}_1; E_2(\gamma, \bar{\delta}, \tilde{\delta}, \varepsilon, x) \bigg\}.$$

This probability is lower bounded by, for any $L$ [condition on $N_i(\varepsilon x)$],

$$\sum_{0 \leq n_i \leq Lx, i \in \mathcal{J}_1} \mathbb{P}\{E_2(\gamma, \bar{\delta}, \tilde{\delta}, \varepsilon, x) \mid N_i(\varepsilon x) = n_i, i \in \mathcal{J}_1\}$$

$$\times \prod_{i \in \mathcal{J}_1} \mathbb{P}\{N_i(\varepsilon x) = n_i\}$$

(B.1)
$$\times \mathbb{P}\bigg\{ (r - c)\bar{\mathbf{A}}_i^r(\varepsilon x)$$

$$- \sum_{k \in \mathcal{J}_1} r_k N_k(\varepsilon x)[\mathbb{E}\{\mathbf{U}_k\} + \bar{\delta}_k] > (1 + \gamma)x, i \in \mathcal{J}_1;$$

$$(r - c)\bar{\mathbf{A}}_i(\varepsilon x) + (r - c)N_i(\varepsilon x)[\mathbb{E}\{\mathbf{A}_i\} - \tilde{\delta}_i]$$

$$- d_i N_i(\varepsilon x)[\mathbb{E}\{\mathbf{U}_i\} + \tilde{\delta}_i]$$

$$- \sum_{k \in \mathcal{J}_1 \setminus \{i\}} r_k N_k(\varepsilon x)[\mathbb{E}\{\mathbf{U}_k\} + \bar{\delta}_k] > (1 + \gamma)x,$$

$$i \in \mathcal{J}_1 \mid N_i(\varepsilon x) = n_i, i \in \mathcal{J}_1 \bigg\}.$$

Before proceeding, we first state a useful lemma (a proof is given at the end of this section). □

LEMMA B.1. *For all $\varepsilon$, $\gamma, \bar{\delta}, \tilde{\delta} > 0$,*

(B.2) $$\mathbb{P}\{E_2(\gamma, \bar{\delta}, \tilde{\delta}, \varepsilon, x) \mid N_i(\varepsilon x) = n_i, i \in \mathcal{J}_1\} \to 1,$$

*as $x \to \infty$ uniformly in $n_i \geq 0, i \in \mathcal{J}_1$, and*

$$\text{(B.3)} \qquad \frac{\mathbb{P}\{N_i(\varepsilon x) = n_i\}}{\mathbb{P}\{\mathbf{A}_i > \varepsilon x\}} \to 1$$

*for all $i \in \mathcal{J}_1$ as $x \to \infty$ uniformly in $0 \leq n_i \leq Lx$.*

Equations (B.2) and (B.3) imply that, for any $L < \infty$ and $\eta > 0$, one can lower bound (B.1) for $x$ large enough by

$$
(1-\eta) \sum_{n_i \leq Lx, i \in \mathcal{J}_1} \mathbb{P}_{\mathcal{J}_0} \Bigg\{ (r-c)\bar{\mathbf{A}}_i^r(\varepsilon x)
$$
$$
- \sum_{k \in \mathcal{J}_1} r_k n_k [\mathbb{E}\{\mathbf{U}_k\} + \delta_k]n > (1+\gamma)x, i \in \mathcal{J}_0;
$$
$$
(r-c)\bar{\mathbf{A}}_i(\varepsilon x) + (r-c)n_i[\mathbb{E}\{\mathbf{A}_i\} - \tilde{\delta}_i]
$$
$$
- d_i n_i [\mathbb{E}\{\mathbf{U}_i\} + \tilde{\delta}_i]
$$
$$
- \sum_{k \in \mathcal{J}_1 \setminus \{i\}} r_k n_k [\mathbb{E}\{\mathbf{U}_k\} + \bar{\delta}_k] > (1+\gamma)x,
$$
$$
i \in \mathcal{J}_1 \mid N_i(\varepsilon x) = n_i, i \in \mathcal{J}_1 \Bigg\}
$$
$$
\times \prod_{i \in \mathcal{J}_1} \mathbb{P}\{\mathbf{A}_i > \varepsilon x\}.
$$

As before, deconditioning upon $\bar{\mathbf{A}}_i$ and $\bar{\mathbf{A}}_i^r$ and applying a similar change of variables as in Appendix A, we obtain the lower bound

$$
(1-\eta) \Bigg( \prod_{i \in \mathcal{J}_1} \frac{p_i}{\mathbb{E}\{\mathbf{A}_i\}} \Bigg)
$$
$$
\times \int_{1 \leq y_i \leq Lx, i \in \mathcal{J}_1} \prod_{i \in \mathcal{J}_0} \mathbb{P} \Bigg\{ (r-c)\mathbf{A}_i^r > (1+\gamma)x + \sum_{k \in \mathcal{J}_1} (r_k - \rho_k + \bar{\delta}_k)y_k \Bigg\}
$$
$$
\times \prod_{i \in \mathcal{J}_1} \mathbb{P} \Bigg\{ (r-c)\mathbf{A}_i > (1+\gamma)x + (d_i - \rho_i + \tilde{\delta}_i)y_i
$$
$$
+ \sum_{k \in \mathcal{J}_1 \setminus \{i\}} (r_k - \rho_k + \bar{\delta}_k)y_k \Bigg\} dy.
$$

Now write

$$
(1-\eta) \int_{1 \leq y_i \leq Lx, i \in \mathcal{J}_1} \cdots = (1-\eta) \int_{y_i \geq 0, i \in \mathcal{J}_1} \cdots - (1-\eta) \int_{\{1 \leq y_i \leq Lx, i \in \mathcal{J}_1\}^c} \cdots
$$

(the complement taken with respect to the nonnegative orthant). The first term in the right-hand side can be handled as in the proof of the upper bound (the only difference is the factor $1 + \gamma$ instead of $1 - \gamma - \delta$). The next lemma shows that the second term can be neglected.

Lemma B.2.

$$
\lim_{L \to \infty} \limsup_{x \to \infty} \frac{1}{P(x)} \int_{\{1 \leq y_i \leq Lx, i \in \mathcal{J}_1\}^c} \prod_{i \in \mathcal{J}_0} \mathbb{P} \bigg\{ (r - c) \mathbf{A}_i^r > (1 + \gamma) x
$$
$$
+ \sum_{k \in \mathcal{J}_1} (r_k - \rho_k + \bar{\delta}_k) y_k \bigg\}
$$
$$
\times \prod_{i \in \mathcal{J}_1} \mathbb{P} \bigg\{ (r - c) \mathbf{A}_i > (1 + \gamma) x + (d_i - \rho_i + \tilde{\delta}_i) y_i
$$
$$
+ \sum_{k \in \mathcal{J}_1 \setminus \{i\}} (r_k - \rho_k + \bar{\delta}_k) y_k \bigg\} \, dy
$$
$$
= 0.
$$

Proof. The integral over the regions in which at least one $y_i$ is smaller than 1 is easily shown to be of $o(P(x))$, so we concentrate on the set $\{0 \leq y_i \leq Lx,$
$i \in \mathcal{J}_1\}^c$. The integral

$$
\int_{\{0 \leq y_i \leq Lx, i \in \mathcal{J}_1\}^c} \prod_{i \in \mathcal{J}_0} \mathbb{P} \bigg\{ (r - c) \mathbf{A}_i^r > (1 + \gamma) x + \sum_{k \in \mathcal{J}_1} (r_k - \rho_k + \bar{\delta}_k) y_k \bigg\}
$$
$$
\times \prod_{i \in \mathcal{J}_1} \mathbb{P} \bigg\{ (r - c) \mathbf{A}_i > (1 + \gamma) x + (d_i - \rho_i + \tilde{\delta}_i) y_i
$$
$$
+ \sum_{k \in \mathcal{J}_1 \setminus \{i\}} (r_k - \rho_k + \bar{\delta}_k) y_k \bigg\} \, dy
$$

is bounded from above by

$$
\bigg( \prod_{i \in \mathcal{J}_0} \mathbb{P} \{ (r - c) \mathbf{A}_i^r > (1 + \gamma) x \}
$$
$$
\times \sum_{j \in \mathcal{J}_1} \int_{y_j \geq Lx, y_i \geq 0, i \in \mathcal{J}_1, i \neq j} \prod_{i \in \mathcal{J}_1} \mathbb{P} \bigg\{ (r - c) \mathbf{A}_i
$$
$$
> (1 + \gamma) x + (d_i - \rho_i + \tilde{\delta}_i) y_i
$$

$$+ \sum_{k \in \mathcal{J}_1 \setminus \{i\}} (r_k - \rho_k + \bar{\delta}_k) y_k, i \in \mathcal{J}_1 \Bigg\} \, dy.$$

Observing that the integrals can be separated, we obtain the upper bound

$$O\Bigg(\prod_{i \in \mathcal{J}_0} \mathbb{P}\{\mathbf{A}_i^r > x\}\Bigg) \sum_{j \in \mathcal{J}_1} O(\mathbb{P}\{\mathbf{A}_j^r > Lx\}) \prod_{i \in \mathcal{J}_1, i \neq j} O\Bigg(\prod_{i \in \mathcal{J}_0} \mathbb{P}\{\mathbf{A}_i^r > x\}\Bigg)$$

$$= O(P(x)) \sum_{j \in \mathcal{J}_1} \frac{\mathbb{P}\{\mathbf{A}_j^r > Lx\}}{\mathbb{P}\{\mathbf{A}_j^r > x\}}.$$

The result then follows immediately. $\quad\square$

PROOF OF LEMMA B.1. Equation (B.2) follows immediately from the following result. Let $\mathbf{S}_n = \mathbf{X}_1 + \cdots + \mathbf{X}_n$ be a random walk with i.i.d. step sizes with $\mathbb{E}\{\mathbf{X}_1\} < 0$. Then

$$\limsup_{x \to \infty} \sup_{n \geq 1} \mathbb{P}\{\mathbf{S}_n > x\} \leq \lim_{x \to \infty} \mathbb{P}\Big\{\sup_{n \geq 1} \mathbf{S}_n > x\Big\} = 0,$$

since $\sup_{n \geq 1} \mathbf{S}_n$ is a proper random variable. Apply this result with $\mathbf{X}_j = \mathbf{U}_{ij} - \mathbb{E}\{\mathbf{U}_i\} - \min\{\bar{\delta}_i, \tilde{\delta}_i\}$ and $\mathbf{X}_j = \mathbb{E}\{\mathbf{A}_i\} - \mathbf{A}_{ij}(\varepsilon x) - \min\{\bar{\delta}_i, \tilde{\delta}_i\}$.

In order to prove (B.3), note that, for $n_i \leq Lx$,

$$\frac{\mathbb{P}\{N_i(\varepsilon x) = n_i\}}{\mathbb{P}\{\mathbf{A}_i > \varepsilon x\}} = \mathbb{P}\{\mathbf{A}_i \leq \varepsilon x\}^{n_i} \leq \mathbb{P}\{\mathbf{A}_i \leq \varepsilon x\}^{Lx} = \left(1 - \frac{o(1)}{x}\right)^{Lx} \to 1,$$

as $x \to \infty$. The last equality holds because $\mathbf{A}_i$ has finite mean. $\quad\square$

**Acknowledgments.** The authors would like to thank Onno Boxma and Miranda van Uitert for useful comments on an earlier version of the paper.

B. ZWART
DEPARTMENT OF MATHEMATICS
AND COMPUTING SCIENCE
EINDHOVEN UNIVERSITY
OF TECHNOLOGY
P.O. BOX 513
5600 MB EINDHOVEN
THE NETHERLANDS
E-MAIL: zwart@win.tue.nl

S. BORST
CWI
P.O. BOX 94079
1090 GB AMSTERDAM
THE NETHERLANDS
E-MAIL: sem@cwi.nl

M. Mandjes
CWI
P.O. Box 94079
1090 GB Amsterdam
The Netherlands
e-mail: michel@cwi.nl